\documentclass[11pt,amsfonts]{article}
\usepackage{graphicx}
\usepackage{latexsym}
\usepackage{amssymb}
\usepackage{amsmath}
\usepackage{layout}
\newtheorem{prop}{Proposition}
\newtheorem{lemma}{Lemma}
\newtheorem{definition}{Definition}

\newtheorem{theorem}{Theorem}
\newtheorem{remark}{Remark}

\def\real{{\mathord{{\rm I\kern-2.8pt R}}}}        
\def\inte{{\mathord{{\rm I\kern-2.8pt N}}}}

\def\sZZ{{\rm Z\kern-2.8ptem{}Z}}

\def\z{{\mathchoice
  {\sZZ}
  {\sZZ}
  {\rm Z\kern-0.30em{}Z}
  {\rm Z\kern-0.25em{}Z} }}
\def\sQQ{{\kern 0.27em \vrule height1.45ex width0.03em depth0em
          \kern-0.30em \rm Q}}
\def\qu{{\mathchoice
    {\sQQ}
    {\sQQ}
  {\kern 0.225em \vrule height1.05ex width0.025em depth0em \kern-0.25em \rm Q}
  {\kern 0.180em \vrule height0.78ex width0.020em depth0em \kern-0.20em \rm Q}
        }}
\def\sCC{{\kern 0.27em \vrule height1.45ex width0.03em depth0em
          \kern-0.30em \rm C}}
\def\complex{{\mathchoice
    {\sCC}
    {\sCC}
  {\kern 0.225em \vrule height1.05ex width0.025em depth0em \kern-0.25em \rm C}
  {\kern 0.180em \vrule height0.78ex width0.020em depth0em \kern-0.20em \rm C}
        }}

\newcommand{\ba}{\begin{array}}
\newcommand{\ea}{\end{array}}
\newcommand{\be}{\begin{equation}}
\newcommand{\ee}{\end{equation}}
\newcommand{\bea}{\begin{eqnarray}}
\newcommand{\eea}{\end{eqnarray}}
\newcommand{\beaa}{\begin{eqnarray*}}
\newcommand{\eeaa}{\end{eqnarray*}}

\newcommand{\eps}{\varepsilon}

\def\b{\beta}

\def\z{\zeta}

\font\tenmath=msbm10 \font\sevenmath=msbm7 \font\fivemath=msbm5
\newfam\mathfam \textfont\mathfam=\tenmath
\scriptfont\mathfam=\sevenmath \scriptscriptfont\mathfam=\fivemath

\def \b{\noindent}

\def \={{\buildrel {\rm (law)} \over =}}

\def\qed{ \hfill \vrule width.25cm height.25cm depth0cm\smallskip}

\newcommand{\basa}{\begin{assumption}}
\newcommand{\easa}{\end{assumption}}

\newcommand{\bas}{\begin{assum}}
\newcommand{\eas}{\end{assum}}

\newcommand{\ignore}[1]{}
\begin{document}

\renewcommand{\thefootnote}{\fnsymbol{footnote}}

\date{ }
\title{Brownian and fractional  Brownian stochastic currents via Malliavin calculus }
\author{Franco Flandoli  $^{1}\quad$Ciprian A. Tudor
$^{2}\vspace*{0.1in}$\\$^{1}$
 Dipartimento di Matematica Applicata, Universita di Pisa,\\ Via Bonnano 25B, I-56126, Pisa, Italy
\\flandoli@dma.unipi.it\vspace*{0.1in}\\$^{2}$SAMOS/MATISSE,
Centre d'Economie de La Sorbonne,\\ Universit\'e de
Panth\'eon-Sorbonne Paris 1,\\90, rue de Tolbiac, 75634 Paris Cedex
13, France.\\tudor@univ-paris1.fr\vspace*{0.1in}}  \maketitle

\maketitle

\begin{abstract}
By using Malliavin calculus and multiple Wiener-It\^o integrals, we study the existence and the regularity of stochastic currents defined as Skorohod (divergence) integrals with respect to the Brownian motion and to the fractional Brownian motion. We consider also the multidimensional multiparameter case and we compare the regularity of the current as a distribution in negative Sobolev spaces with its regularity in Watanabe space.

\end{abstract}

\vskip0.5cm

{\bf  2000 AMS Classification Numbers:} 60G15, 60G18,  60H05, 76M35, 60H05

 \vskip0.3cm


\section{Introduction}

The concept of current is proper of geometric measure theory. The simplest
example is the functional
\[
\varphi\mapsto\int_{0}^{T}\left\langle \varphi\left(  \gamma\left(  t\right)
\right)  ,\gamma^{\prime}\left(  t\right)  \right\rangle _{\mathbb{R}^{d}}dt
\]
defined over the set of all smooth compact support vector fields
$\varphi:\mathbb{R}^{d}\rightarrow\mathbb{R}^{d}$, with $\gamma:\left[
0,T\right]  \rightarrow$ $\mathbb{R}^{d}$ being a rectifiable curve. This
functional defines a vector valued distribution. Let us denote it by
\begin{equation}
\int_{0}^{T}\delta\left(  x-\gamma\left(  t\right)  \right)  \gamma^{\prime
}\left(  t\right)  dt.\label{deterministic current}%
\end{equation}
We address for instance to the books \cite{Fed}, \cite{Simon}, \cite{Morgan},
\cite{GMS}, for definitions, results and applications.

The stochastic analog of 1-currents is a natural concept, where the
deterministic curve $\left(  \gamma\left(  t\right)  \right)  _{t\in\left[
0,T\right]  }$ is replaced by a stochastic process $\left(  X_{t}\right)
_{t\in\left[  0,T\right]  }$ in $\mathbb{R}^{d}$ and the stochastic integral
must be properly interpreted. Several works deal with random currents, see for
instance \cite{Manabe1}, \cite{Ikeda}, \cite{Ochi}, \cite{IO}, \cite{MO},
\cite{Manabe2}, \cite{Kuwada1}  \cite{Gia}, \cite{Kuwada2}, \cite{FGR}.

The difference between classical integration theory and random currents is the
attempt to understand the latter as random distribution in the strong sense:
random variables taking values in the space of distributions. The question is
not simply how to define $\int_{0}^{1}\left\langle \varphi\left(
X_{t}\right)  ,dX_{t}\right\rangle _{\mathbb{R}^{d}}$ for every given test
functions, but when this operation defines, \textit{for almost every
realization} of the process $X$, a continuous functional on some space of test
functions. This is related to (a particular aspect of) T. Lyons theory of
rough paths:\ a random current of the previous form is a concept of pathwise stochastic
integration. The degree of regularity of the random analog of expression
(\ref{deterministic current}) is a fundamental issue.

In this work we study  the mapping defined as
\begin{equation}
\label{xi} \xi (x) = \int _{[0,T]^{N}} \delta (x-B_{s} )  dB_{s},
\hskip0.5cm x\in \mathbb{R} ^{d}, T>0.
\end{equation}
Here the integrator process is a $d$- dimensional Wiener process with multidimensional
time parameter or a $d$-dimensional fractional Brownian motion.

A first direction of investigation is the regularity with respect to
the  variable $x  \in \mathbb{R}^{d} $ of the mapping  given by
(\ref{xi})  in the (deterministic) Sobolev spaces $H^{-r}
(\mathbb{R}^{d}; \mathbb{R}^{d})$.
This is the dual space to $H^{r}\left(  \mathbb{R}^{d},\mathbb{R}^{d}\right)
$, or equivalently the space of all vector valued distributions $\varphi$ such that
$\int_{\mathbb{R}^{d}}\left(  1+\left|  x\right|  ^{2}\right)  ^{-r}\left|
\widehat{\varphi}\left(  x\right)  \right|  ^{2}dx<\infty$, $\widehat{\varphi
}$ being the Fourier transform of $\varphi$.  Some results on this have been
recently obtained by \cite{Gia}, \cite{FG} or \cite{FGR} using
different techniques of the stochastic calculus. We propose here a
new approach to this problem. The main tool is constituted by the
Malliavin calculus based on the {\em Wiener-It\^o  chaos expansion.
} Actually the delta Dirac function $\delta (x-B_{s} )$ can be
understood as a distribution in the Watanabe sense (see \cite{Wa})
and it has been the object of study by several authors, as
\cite{NV}, \cite{Imk1} or \cite{Imk2}  in the Brownian motion  case,
or \cite{CNT},  \cite{Ed} or \cite{Hu} in the fractional Brownian
motion case. In particular, it is possible to obtain the
decomposition of the delta Dirac function into an orthogonal sum of  multiple
Wiener-It\^o integrals and as a consequence it is easy to get the
chaos expansion of the integral (\ref{xi}) since the divergence
integral acts as a creation operator on Fock spaces (see Section 2).
As it will be seen, the chaos expansion obtained will be useful to
obtain the regularity properties of the mapping $\xi$. Another
advantage of this method is the fact that it can be relatively
easily extended to multidimensional settings.

Besides the estimation of the Sobolev regularity with respect to
$x$, we are also interested to study the regularity with respect
to $\omega$, that is, as a functional in the Watanabe sense, of the
integral (\ref{xi}).  Our interest comes from the following
observation. Consider $N=d=1$.  It is well-known (see Nualart and
Vives \cite{NV} ) that for fixed $x\in \mathbb{R}$ we have
\begin{equation*}
\delta (x-B_{s}) \in \mathbb{D}^{-\alpha , 2} \hskip0.5cm \mbox{for
any } \alpha >\frac{1}{2}
\end{equation*}
where $\mathbb{D}^{-\alpha , 2}$ are the Watanabe (or
Sobolev-Watanabe) spaces introduced in Section 2. On the other side, more or less surprisedly   the same order of
regularity holds with respect to $x$; for fixed $\omega$, the
mapping $g(x)=\delta (x-B_{s}(\omega ))$ belongs to the negative
Sobolev space $H^{-r} (\mathbb{R}; \mathbb{R}) $  for every
$r>\frac{1}{2}$. Indeed,  since the Fourier transform of $g$ is
$\hat{g} (x)= e^{-ixB_{s}} $ we have
\begin{equation*}
\left| g\right| ^{2} _{H^{-r}(\mathbb{R}; \mathbb{R})} = \int
_{\mathbb{R}} \vert \hat{g}(x) \vert ^{2} (1+ x^{2})^{-r}dx
\end{equation*}
and this is finite if and only if $r>\frac{1}{2}$.

One can ask the question if this similarity of the order of
regularity  in the deterministic and stochastic Sobolev space still
holds  for the functional $\xi$ defined by \ref{xi}; and actually in the case of dimension $d=1$ the above  property still holds: we have the same regularity of $\xi$ both with respect to $x$ and with respect to $\omega$.

 One can moreover ask if it holds  for other Gaussian processes. The answer to this question is negative,
since we show that (actually, this has also been proved in \cite{FGR} but using another integral with respect to fBm) in the fractional Brownian motion case the
mapping (\ref{xi}) belongs to $H^{-r} (\mathbb{R}; \mathbb{R}) $ for any $r>\frac{1}{2H}-\frac{1}{2} $ (as a
function in $x$) and to $\mathbb{D}^{-\alpha ,2} $ with $\alpha
>\frac{3}{2}-\frac{1}{2H}$ (as a function on $\omega$).

We organized our paper as follows. Section 2 contains some preliminaries on Malliavin calculus and multiple Wiener-It\^o integrals. Section 3 contains a discussion about random distribution where we unify the definition  of the quantity $\delta (x-X_{s}(\omega))$ ($X$ is a Gaussian process on $\mathbb{R}^{d}$) which in principle can be understood as a distribution with respect to $x$ and also as a distribution in the Watanabe sense when it is regarded as o function of $\omega$. In Section 4 we study the existence and the regularity of the stochastic currents driven by  a $N$ parameter Brownian motion in $\mathbb{R}^{d}$ while in Section 5 concerning the same problem when the driving process is the fractional Brownian motion. In Section 6 we give the regularity of the integral (\ref{xi}) with respect to $\omega$ and we compare it with its regularity in $x$.
\vskip0.5cm

\section{Preliminaries}

Here we describe the elements from stochastic analysis that we will
need in the paper. Consider ${\mathcal{H}}$ a real separable Hilbert
space and $(B (\varphi), \varphi\in{\mathcal{H}})$ an isonormal
Gaussian process on a probability space $(\Omega, {\cal{A}}, P)$, that is a centered Gaussian family of random
variables such that $\mathbf{E}\left( B(\varphi) B(\psi) \right)  =
\langle\varphi, \psi\rangle_{{\mathcal{H}}}$.
Denote by $I_{n}$ the multiple stochastic integral with respect to
$B$ (see \cite{N}). This $I_{n}$ is actually an isometry between the
Hilbert space ${\mathcal{H}}^{\odot n}$(symmetric tensor product)
equipped with the scaled norm
$\frac{1}{\sqrt{n!}}\Vert\cdot\Vert_{{\mathcal{H}}^{\otimes n}}$ and
the Wiener chaos of order $n$ which is defined as the closed linear
span of the random variables $H_{n}(B(\varphi))$ where
$\varphi\in{\mathcal{H}}, \Vert\varphi\Vert_{{\mathcal{H}}}=1$ and
$H_{n}$ is the Hermite polynomial of degree $n\geq 1$
\begin{equation*}
H_{n}(x)=\frac{(-1)^{n}}{n!} \exp \left( \frac{x^{2}}{2} \right)
\frac{d^{n}}{dx^{n}}\left( \exp \left( -\frac{x^{2}}{2}\right)
\right), \hskip0.5cm x\in \mathbb{R}.
\end{equation*}

 The isometry of multiple integrals can be written as: for $m,n$ positive integers,
\begin{eqnarray}
\mathbf{E}\left(I_{n}(f) I_{m}(g) \right) &=& n! \langle f,g\rangle _{{\mathcal{H}}^{\otimes n}}\quad \mbox{if } m=n,\nonumber \\
\mathbf{E}\left(I_{n}(f) I_{m}(g) \right) &= & 0\quad \mbox{if } m\not=n.\label{iso}
\end{eqnarray}
It also holds that
\begin{equation*}
I_{n}(f) = I_{n}\big( \tilde{f}\big)
\end{equation*}
where $\tilde{f} $ denotes the symmetrization of $f$ defined by $\tilde{f}%
(x_{1}, \ldots , x_{x}) =\frac{1}{n!} \sum_{\sigma \in {\cal S}_{n}}
f(x_{\sigma (1) }, \ldots , x_{\sigma (n) } ) $.

We recall that any square integrable random variable which is
measurable with respect to the $\sigma$-algebra generated by $B$ can
be expanded into an orthogonal sum of multiple stochastic integrals
\begin{equation}
\label{sum1} F=\sum_{n\geq0}I_{n}(f_{n})
\end{equation}
where $f_{n}\in{\mathcal{H}}^{\odot n}$ are (uniquely determined)
symmetric functions and $I_{0}(f_{0})=\mathbf{E}\left[  F\right]  $.

Let $L$ be the Ornstein-Uhlenbeck operator
\begin{equation*}
LF=-\sum_{n\geq 0} nI_{n}(f_{n})
\end{equation*}
if $F$ is given by (\ref{sum1}).

For $p>1$ and $\alpha \in \mathbb{R}$ we introduce the
Sobolev-Watanabe space $\mathbb{D}^{\alpha ,p }$  as the closure of
the set of polynomial random variables with respect to the norm
\begin{equation*}
\Vert F\Vert _{\alpha , p} =\Vert (I -L) ^{\frac{\alpha }{2}} \Vert
_{L^{p} (\Omega )}
\end{equation*}
where $I$ represents the identity. In this way, a random variable
$F$  as in (\ref{sum1}) belongs $\mathbb{D}^{\alpha , 2}$ if and
only if
\begin{equation*}
\sum_{n\geq 0} (1+n) ^{\alpha } \Vert I_{n}(f_{n}) \Vert
_{L^{2}(\Omega)} ^{2} =\sum_{n\geq 0} (1+n) ^{\alpha }n! \Vert f_{n}
\Vert ^{2} _{{\cal{H}}^{\otimes n}}<\infty .
\end{equation*}
We denote by $D$  the Malliavin  derivative operator that acts on smooth functions of the form $F=g(B(\varphi _{1}), \ldots , B(\varphi_{n}))$ ($g$ is a smooth function with compact support and $\varphi_{i} \in {{\cal{H}}}$)
\begin{equation*}
DF=\sum_{i=1}^{n}\frac{\partial g}{\partial x_{i}}(B(\varphi _{1}), \ldots , B(\varphi_{n}))\varphi_{i}.
\end{equation*}
The operator $D$ is continuous from $\mathbb{D} ^{\alpha , p} $ into
$\mathbb{D} ^{\alpha -1, p} \left( {\cal{H}}\right).$ The adjoint of
$D$ is denoted by $\delta $ and is called the divergence (or
Skorohod) integral. It is a continuous operator from $\mathbb{D}
^{\alpha, p } \left( {\cal{H}}\right)$ into $\mathbb{D} ^{\alpha -1,
p}$. For adapted integrands, the divergence integral coincides to
the classical It\^o integral. We will use the notation
\begin{equation*}
\delta (u) =\int_{0}^{T} u_{s} dB_{s}.
\end{equation*}
Let  $u$ be  a stochastic process having the chaotic decomposition $u_{s}=\sum _{n\geq 0} I_{n}(f_{n}(\cdot ,s))$
where $f_{n}(\cdot, s)\in {\cal{H}}^{\otimes n}$ for every $s$. One
can prove that $u \in {\rm Dom} \ \delta$ if and only if $\tilde f_n
\in {\cal{H}}^{\otimes (n+1)}$ for every $n \geq 0$, and
$\sum_{n=0}^{\infty}I_{n+1}(\tilde f_n)$ converges in $L^2(\Omega)$.
In this case,
$$\delta(u)=\sum_{n=0}^{\infty}I_{n+1}(\tilde f_n) \quad \mbox{and}
\quad \mathbf{E}|\delta(u)|^2=\sum_{n=0}^{\infty}(n+1)! \ \|\tilde
f_n\|_{{\cal{H}}^{\otimes (n+1)}}^{2}.$$

In the present work we will consider divergence integral with respect to a Brownian motion in $\mathbb{R}^{d}$ as well as with respect to a fractional Brownian motion in $\mathbb{R}^{d}$.

Throughout this paper we will denote  by $p_{s}(x) $ the Gaussian
kernel of variance $s>0$ given by $p_{s}(x)= \frac{1}{\sqrt{2\pi s}} e^{-\frac{x^{2}}{2s}},
x\in \mathbb{R}$ and for $x=(x_{1}, \ldots , x_{d})\in \mathbb{R}^{d}$ by $p_{s} ^{d}(x)= \prod_{i=1}^{d}p_{s}(x_{i}).$

\section{Random distributions}

We study now the regularity of the stochastic integral given by
(\ref{xi}). Our method is based on the Wiener-It\^o chaos
decomposition. Let $X$ be isonormal Gaussian process with variance
$R(s,t)$. We will use the following decomposition of the delta Dirac
function (see Nualart and Vives \cite{NV}, Imkeller et al. \cite{Imk1}, Eddahbi et al. \cite{Ed}) into orthogonal multiple Wiener-It\^o
integrals
\begin{equation}
\label{deltagen} \delta (x-X_{s}) = \sum _{n\geq 0} R(s)
^{-\frac{n}{2}}p_{R(s)}(x) H_{n} \left(
\frac{x}{R(s)^{\frac{1}{2}}}\right) I_{n}\left( 1_{[0,s]}^{\otimes
n}\right)
\end{equation}
where $R(s):=R(s,s)$, $p_{R(s)}$ is the Gaussian kernel of variance
$R(s)$, $H_{n}$ is the Hermite polynomial of degree $n$ and $I_{n}$
represents the multiple Wiener-It\^o integral of degree $n$ with
respect to the Gaussian process $X$ as defined in the previous section.

A key element of the entire work is the quantity $\delta (x-B_{s})$. This can be understood as a generalized random variable in some Sobolev-Watanbe space as well as a  generalized function with respect to the variable $x$. The purpose of this section is to give an unitary definition of  the delta Dirac as a random distribution.

Let $D\left(  \mathbb{R}^{d}\right)  $ be the space of smooth compact support
functions on $\mathbb{R}^{d}$ and let $D^{\prime}\left(  \mathbb{R}%
^{d}\right)  $ be its dual, the space of distributions, endowed with the usual
topologies. We denote by $\left\langle S,\varphi\right\rangle $ the dual
pairing between $S\in D^{\prime}\left(  \mathbb{R}^{d}\right)  $ and
$\varphi\in D\left(  \mathbb{R}^{d}\right)  $. We say that a distribution
$S\in D^{\prime}\left(  \mathbb{R}^{d}\right)  $ is of class $L_{loc}%
^{p}\left(  \mathbb{R}^{d}\right)  $, $p\geq1$, if there is $f\in L_{loc}%
^{1}\left(  \mathbb{R}^{d}\right)  $  such that $\left\langle S,\varphi
\right\rangle =\int_{\mathbb{R}^{d}}f\left(  x\right)  \varphi\left(
x\right)  dx$ for every $\varphi\in D\left(  \mathbb{R}^{d}\right)  $. In this
case we also say that the distribution is given by a function.

Let $\left(  \Omega,\mathcal{A},P\right)  $ be a probability space, with
expectation denoted by $\mathbf{E}.$

\begin{definition}
We call random distribution (on $\mathbb{R}^{d}$, based on $\left(
\Omega,\mathcal{A},P\right)  $) a measurable mapping $\omega\mapsto S\left(
\omega\right)  $ from $\left(  \Omega,\mathcal{A}\right)  $ to $D^{\prime
}\left(  \mathbb{R}^{d}\right)  $ with the Borel $\sigma$-algebra.
\end{definition}

Given a random distribution $S$, for every $\varphi\in D\left(  \mathbb{R}%
^{d}\right)  $ the real valued function $\omega\mapsto\left\langle S\left(
\omega\right)  ,\varphi\right\rangle $ is measurable. The converse is also true.

Similarly to the deterministic case, we could say that a random distribution
$S$ is of class $L^{0}\left(  \Omega,L_{loc}^{1}\left(  \mathbb{R}^{d}\right)
\right)  $ if there exists a measurable function $f:\Omega\times\mathbb{R}%
^{d}\rightarrow\mathbb{R}$, $f\left(  \omega,\cdot\right)  \in L_{loc}%
^{1}\left(  \mathbb{R}^{d}\right)  $ for $P$-a.e. $\omega\in\Omega$, such that
$\left\langle S\left(  \omega\right)  ,\varphi\right\rangle =\int
_{\mathbb{R}^{d}}f\left(  \omega,x\right)  \varphi\left(  x\right)  dx$ for
every $\varphi\in D\left(  \mathbb{R}^{d}\right)  $. In this case we also say
that the distribution is given by a random field. However, this concept is
restrictive if we have to deal with true distributions.

There is an intermediate concept, made possible by the simultaneous presence
of the two variables $\omega\in\Omega$ and $x\in\mathbb{R}^{d}$. One could
have that the random distribution is given by a function with respect to the
$x$ variable, but at the price that it is distribution-valued (in the
Sobolev-Watanabe sense)\ in the $\omega$ variable. To this purpose, assume
that a real separable Hilbert space $\mathcal{H}$ is given, an isonormal
Gaussian process $\left(  W\left(  h\right)  ,h\in\mathcal{H}\right)  $ on
$\left(  \Omega,\mathcal{A},P\right)  $ is given. For $p>1$ and $\alpha
\in\mathbb{R}$, denote by $\mathbb{D}^{\alpha,p}$ the Sobolev-Watanabe space
of generalized random variables on $\left(  \Omega,\mathcal{A},P\right)  $,
defined in terms of Wiener chaos expansion.

\begin{definition}
Given $p\geq1$ and $\alpha\leq0$, a random distribution $S$ is of class
$L_{loc}^{1}\left(  \mathbb{R}^{d};\mathbb{D}^{\alpha,p}\right)  $ if
$\mathbf{E}\left[  \left|  \left\langle S,\varphi\right\rangle \right|
^{p}\right]  <\infty$ for all $\varphi\in D\left(  \mathbb{R}^{d}\right)  $
and there exists a function $f\in L_{loc}^{1}\left(  \mathbb{R}^{d}%
;\mathbb{D}^{\alpha,p}\right)  $ such that
\[
\left\langle S,\varphi\right\rangle =\int_{\mathbb{R}^{d}}\varphi\left(
x\right)  f\left(  x\right)  dx
\]
for every $\varphi\in D\left(  \mathbb{R}^{d}\right)  $.
\end{definition}

Let us explain the definition. The integral $\int_{\mathbb{R}^{d}}%
\varphi\left(  x\right)  f\left(  x\right)  dx$ is of Bochner type (the
integral of the $\mathbb{D}^{\alpha,p}$-valued function $x\mapsto
\varphi\left(  x\right)  f\left(  x\right)  $). The integral $\int
_{\mathbb{R}^{d}}\varphi\left(  x\right)  f\left(  x\right)  dx$ is thus, a
priori, an element of $\mathbb{D}^{\alpha,p}$. The random variable
$\omega\mapsto\left\langle S\left(  \omega\right)  ,\varphi\right\rangle $,
due to the assumption $\mathbf{E}\left[  \left|  \left\langle S,\varphi
\right\rangle \right|  ^{p}\right]  <\infty$, is an element of $\mathbb{D}%
^{\alpha^{\prime},p}$ for all $\alpha^{\prime}\leq0$. The definition requires
that they coincide, as elements of $\mathbb{D}^{\alpha,p}$, a priori. A
fortiori, since $\omega\mapsto\left\langle S\left(  \omega\right)
,\varphi\right\rangle $ is an element of $\mathbb{D}^{0,p}$, the same must be
true for $\int_{\mathbb{R}^{d}}\varphi\left(  x\right)  f\left(  x\right)
dx$. Thus, among the consequences of the definition there is the fact that
$\int_{\mathbb{R}^{d}}\varphi\left(  x\right)  f\left(  x\right)
dx\in\mathbb{D}^{0,p}$ although $f\left(  x\right)  $ lives only in
$\mathbb{D}^{\alpha,p}$.

The following theorem treats our main example. Given any $d$- dimensional random variable $X$
on $\left(  \Omega,\mathcal{A},P\right)  $, let $S$ be the random distribution
defined as
\[
\left\langle S\left(  \omega\right)  ,\varphi\right\rangle :=\varphi\left(
X\left(  \omega\right)  \right), \hskip0.5cm \varphi \in D(\mathbb{R}^{d}).
\]
We denote this random distribution by
\[
\delta\left(  x-X\left(  \omega\right)  \right)
\]
and we call it the delta-Dirac at $X\left(  \omega\right)  $. Denote Hermite
polynomials by $H_{n}\left(  x\right)  $, multiple Wiener integrals by $I_{n}%
$, Gaussian kernel of variance $\sigma^{2}$ by $p_{\sigma^{2}}\left(
x\right)  $.

 If ${\mathcal{H}}_{1}, \ldots , {\mathcal{H}}_{d}$ are real and separable Hilbert spaces, a $d$- dimensional  isonormal process is defined on the product space $(\Omega, {\cal{A}}, P)$, $\Omega =\Omega_{1}\times \ldots \times \Omega _{d}, {\cal{A}}={\cal{A}}_{1}\otimes ...\otimes {\cal{A}}_{d},  P=P_{1}\otimes ...\otimes P_{d}$ as a vector with independent components $\left( (W^{1}(h _{1}),\ldots , W^{d}(\varphi _{d})), h_{1}\in \mathcal{H}_{1}, \ldots ,h _{d}\in \mathcal{H}_{d}\right)$ where for every $i=1, \ldots ,d$, $(W^{i}(h_{i}),h _{i}\in {\mathcal{H}}_{i})$ is an  one dimensional isonormal process on $(\Omega _{i}, {\cal{A}}_{i}, P_{i})$. We denote by $\mathbf{E}_{i}$ the expectation on $(\Omega _{i}, {\cal{A}}_{i}, P_{i})$ and by $\mathbf{E}$ the expectation on $(\Omega, {\cal{A}}, P)$.

\begin{theorem}
Let $\mathcal{H}_{i}$, $i=1,\ldots, d$ be real separable Hilbert spaces. Consider  \\
$\left( (W^{1}(h _{1}),\ldots , W^{d}(\varphi _{d})), h_{1}\in \mathcal{H}_{1}, \ldots ,h _{d}\in \mathcal{H}_{d}\right)$ a $d$-dimensional isonormal
Gaussian process on $\left(  \Omega,\mathcal{A},P\right)  $ as above. Then, for every
$h_{i}\in\mathcal{H}_{i}$, $i=1,\ldots , d$,  the random distribution $\delta\left(  \cdot-W\left(
h\right)  \right)  $ is of class $L_{loc}^{1}\left(  \mathbb{R}^{d}%
;\mathbb{D}^{\alpha,2}\right)  $ for some $\alpha<0$ and the associated
element $f$ of $L_{loc}^{1}\left(  \mathbb{R}^{d};\mathbb{D}^{\alpha
,2}\right)  $ is%
\begin{equation*}
f(x)=  \sum_{n_{1},..., n_{d}\geq 0}\prod_{i=1}^{d} \vert h_{i} \vert ^{-n_{i}}p_{\vert h_{i}\vert ^{2}}(x_{i}) H_{n_{i}}\left( \frac{x_{i}}{\vert h_{i}\vert }\right)I^{i}_{n_{i}}\left(  h_{i}^{\otimes n_{i}}\right)
\end{equation*}
where  $x=(x_{1}, \ldots , x_{d}) \in \mathbb{R}^{d}$ and $I^{i}_{n_{i}}$, $i=1,...,d$ denotes the multiple Wiener-It\^o integral with respect to the Wiener process $W^{(i)}$.

We denote this $\mathbb{D}^{\alpha,2}$-valued function $f\left(  x\right)  $
by $\delta\left(  x-W\left(  h\right)  \right)  $.
\end{theorem}
{\bf Proof: } Let us consider first the case $d=1$. In this case
\[
f\left(  x\right)  =\sum_{n\geq0}\left|  h\right|  ^{-n}p_{\left|  h\right|
^{2}}\left(  x\right)  H_{n}\left(  \frac{x}{\left|  h\right|  }\right)
I_{n}\left(  h^{\otimes n}\right)  .
\]
We have to prove that
\[
\varphi\left(  W\left(  h\right)  \right)  =\int_{\mathbb{R}}%
\varphi\left(  x\right)  f\left(  x\right)  dx
\]
for every test function $\varphi$. We have
\begin{eqnarray*}
\int_{\mathbb{R}}\varphi\left(  x\right)  f\left(  x\right)  dx  &=&\sum_{n\geq0}\left|  h\right|  ^{-n}I_{n}\left(  h^{\otimes n}\right)
\int_{\mathbb{R}}\varphi\left(  x\right)  p_{1}\left(  \frac{x}{\left|
h\right|  }\right)  H_{n}\left(  \frac{x}{\left|  h\right|  }\right)  dx\\
& =&\sum_{n\geq0}\vert h\vert ^{-n+1}I_{n}\left(  h^{\otimes n}\right)  \int_{\mathbb{R}}\varphi\left(  \left|  h\right|  x\right)  H_{n}\left(  x\right)
p_{1}\left(  x\right)  dx.
\end{eqnarray*}
On the other hand, by Stroock's formula we can write
\begin{eqnarray*}
\varphi \left( W(h) \right) &=&
\sum_{n\geq 0} \frac{1}{n!} I_{n}\left( D^{(n)}\varphi \left( W(h) \right)\right)=\sum_{n\geq 0} \frac{1}{n!} I_{n}(h^{\otimes n} ) \mathbf{E}\left( \varphi^{(n)}(W(h)) \right)
\end{eqnarray*}
and
\begin{eqnarray*}
\mathbf{E}\left( \varphi^{(n)}(W(h)) \right)&=& \int _{\mathbb{R}} \varphi ^{(n)} (x) p_{\vert h\vert ^{2}}(x) dx
=\vert h\vert \int_{\mathbb{R}} \varphi ^{(n)} (x\vert h\vert ) p_{1}(x) dx \\
&=&(-1) \int_{\mathbb{R}} \varphi ^{(n-1)} (x\vert h\vert ) (p_{1} (x) ) ^{\prime } dx \\
&=& \ldots \\
&=&(-1) ^{n}\vert h\vert ^{n-1} \int_{\mathbb{R}} \varphi (x\vert h\vert) (p_{1} (x) ) ^{(n)} dx.
\end{eqnarray*}
By using the recurrence formula for the Gaussian kernel $p_{1}$
\begin{equation*}
(p_{1} (x) ) ^{(n)}=(-1)^{n} n! p_{1}(x) H_{n}(x)
\end{equation*}
we  get
\begin{equation*}
\varphi\left( W(h) \right)= \sum_{n\geq 0} \int_{\mathbb{R}} \varphi (x\vert h\vert ) H_{n}(x) p_{1}(x) dx=\int_{\mathbb{R}}\varphi(x) f(x)dx.
\end{equation*}
Concerning the case $d\geq 2$, note that
\begin{eqnarray*}
\int_{\mathbb{R}^{d}} \varphi(x)f(x)dx
&=& \sum_{n_{1}, \ldots , n_{d}\geq 0} \prod _{i=1}^{d} \left( \vert h_{i} \vert ^{-n_{i}}I^{i}_{n_{i}}\left(  h_{i}^{\otimes n_{i}}\right)\int_{\mathbb{R}^{d}}dx \varphi (x)p_{\vert h_{i}\vert ^{2}}(x_{i}) H_{n_{i}}\left( \frac{x_{i}}{\vert h_{i}\vert }\right)\right)\\
&=& \sum_{n_{1}, \ldots , n_{d}\geq 0} \prod _{i=1}^{d} \left( \vert h_{i} \vert ^{-n_{i}+1}I^{i}_{n_{i}}\left(  h_{i}^{\otimes n_{i}}\right)\int_{\mathbb{R}^{d}}\varphi(\vert h_{1}\vert x_{1}, \ldots , \vert h_{d}\vert x_{d}) p_{1}(x_{i}) H_{n_{i}}(x_{i}) dx\right).
\end{eqnarray*}
Let us apply Stroock 's formula in several steps. First we apply Stroock formula with respect to the component $W^{i}$
\begin{eqnarray*}
\varphi (W(h))=\varphi (W^{1}(h_{1}), \ldots , W^{d}(h_{d}) ) &=&\sum_{n_{1}\geq 0} \frac{1}{n_{1}!} I^{1}_{n_{1}} \left( \mathbf{E}_{1} \left[ D^{(n), 1}\varphi (W^{1}(h_{1}), \ldots , W^{d}(h_{d}) )\right] \right)\\
&=&\sum_{n_{1}\geq 0} \frac{1}{n_{1}!} I^{1}_{n_{1}}(h_{1}^{\otimes n_{1}}) \mathbf{E}_{1} \left( \frac{\partial^{n_{1}}\varphi }{\partial x_{1}^{n_{1}}} (W^{1}(h_{1}), \ldots , W^{d}(h_{d}) )\right)
\end{eqnarray*}
and then with respect to $W^{2}$
\begin{equation*}
\frac{\partial^{n_{1}}\varphi }{\partial x_{1}^{n_{1}}} (W^{1}(h_{1}), \ldots , W^{d}(h_{d}) )=\sum_{n_{2}\geq 0} \frac{1}{n_{2}!} I^{2}_{n_{2}}(h_{2}^{\otimes n_{2}})\mathbf{E}_{2} \left( \frac{\partial ^{n_{1}+n_{2}}\varphi}{\partial x_{1}^{n_{1}}x_{2}^{n_{2}}}\right) (W^{1}(h_{1}), \ldots , W^{d}(h_{d}) ).
\end{equation*}
We will obtain
\begin{eqnarray*}
\varphi (W(h))&=& \sum_{n_{1}, \ldots , n_{d}\geq 0}\prod_{i=1}^{d} \frac{1}{n_{i}!} I^{i}_{n_{i}}(h_{i}^{\otimes n_{i}})\mathbf{E}\left( \frac{\partial^{n_{1}+..+n_{d}}\varphi}{\partial x_{1}^{n_{1}}..x_{d}^{n_{d}}} (W^{1}(h_{1}), \ldots , W^{d}(h_{d}) )\right)\\
&=&\sum_{n_{1}, \ldots , n_{d}\geq 0}\prod_{i=1}^{d} \frac{1}{n_{i}!} I^{i}_{n_{i}}(h_{i}^{\otimes n_{i}})\int_{\mathbb{R}^{d}}\frac{\partial^{n_{1}+..+n_{d}}\varphi}{\partial x_{1}^{n_{1}}..x_{d}^{n_{d}}}(x)p_{\vert h_{i}\vert ^{2}}(x_{i})dx
\end{eqnarray*}
and now it suffices to follow the case $d=1$ by integrating by parts and using the recurrence formula for the Gaussian kernel.\qed

\section{Regularity of Brownian currents with respect to $x$ }

Let $(B_{t})_{t\in [0,T]^{N}}$ be a Gaussian isonormal process. For every $n\geq 1$, $x\in \mathbb{R}$ and $s=(s_{1}, \ldots ,s_{N})\in [0,T]^{N}$ we will use  the notation

\begin{equation}
\label{an} a_{n}^{x}(s)= R(s) ^{-\frac{n}{2}}p_{R(s)}(x) H_{n}
\left( \frac{x}{R(s)^{\frac{1}{2}}}\right)
\end{equation}
with $R(s)=\mathbf{E}B_{s}^{2}=s_{1}\ldots s_{N}$.
We will start by the following very useful calculation which plays an important role in the sequel.

\begin{lemma}
Let $a_{n}^{x} $ be given by (\ref{an}) and denote by
 $a_{n}^{\hat{x}}(s)$ the Fourier transform of the function $x\to
a_{n}^{x}(s)$. Then it holds  that
\begin{equation}
\label{imp} a_{n}^{\hat{x}}(s)=e^{-\frac{x^{2}R(s)}
{2}}\frac{(-i)^{n}x^{n}}{n!}
\end{equation}
\end{lemma}
{\bf Proof: } Using formula (\ref{deltagen}) one can prove  that that Fourier transform of
$g(x)=\delta (x -B_{s}) $  admits the chaos expansion
\begin{equation}
\label{r1} \hat{g}(x)=\sum_{n\geq 0}a_{n}^{\hat{x}}(s)I_{n}\left(
1_{[0,s]}^{\otimes n}\right)
\end{equation}
where  $a_{n}^{\hat{x}}(s)$ denotes the Fourier transform of the
function $x\to a_{n}^{x}(s)$.

But on the other hand, we know that $\hat{ \delta } (x-B_{s}) =
e^{-ixB_{s}}$  and  using the Stroock's formula to decompose in
chaos a square integrable random variable $F$ infinitely
differentiable in the Malliavin sense
\begin{equation*}
F=\sum _{n\geq 0} \frac{1}{n!} I_{n}\left[ \mathbf{E}(D^{(n)} F)\right]
\end{equation*}
where $D^{(n)} $ denotes the $n$th Malliavin derivative and by the
trivial relation
\begin{equation*}
D_{u}e^{-ixB_{s}}= - ixe^{-ixB_{s}}1_{[0,s]}(u)
\end{equation*}
which implies
\begin{equation*}
D^{(n)} _{u_{1},\ldots , u_{n}}e^{-ixB_{s}}=(-ix)
^{n}e^{-ixB_{s}}1_{[0,s]}^{\otimes n}(u_{1}, \ldots , u_{n})
\end{equation*}
 one obtains
\begin{equation}
\label{r2}
 \hat{g}(x)=\mathbf{E}(e^{-ixB_{s}})\sum_{n\geq 0}\frac{(-i)^{n}x^{n}}{n!}I_{n}\left(
1_{[0,s]}^{\otimes n}\right)= e^{-\frac{x^{2}R(s)} {2}}\sum_{n\geq
0}\frac{(-i)^{n}x^{n}}{n!}I_{n}\left( 1_{[0,s]}^{\otimes n}\right).
\end{equation}
Now, by putting together  (\ref{r1}) and (\ref{r2}) we get  (\ref{imp}). \qed

\subsection{The one-dimensional Brownian currents}

 Let  $B$ be in this part a   the Wiener process. Then by (\ref{an})
\begin{equation}
\label{deltaB} \delta (x-B_{s}) = \sum _{n\geq 0}
s^{-\frac{n}{2}}p_{s}(x) H_{n} \left(
\frac{x}{s^{\frac{1}{2}}}\right) I_{n}\left( 1_{[0,s]}^{\otimes
n}\right)= \sum _{n\geq 0} a_{n}^{x} (s)I_{n}\left(
1_{[0,s]}^{\otimes n}\right),
\end{equation}
here  $I_{n}$ representing the multiple
Wiener-It\^o integral of degree $n$ with respect to the Wiener
process. In this case we have
\begin{equation*}
a_{n}^{x}(s)=s^{-\frac{n}{2}}p_{s}(x) H_{n} \left(
\frac{x}{s^{\frac{1}{2}}}\right).
\end{equation*}

We prove first the following result, which is already known (see
\cite{Gia}, \cite{FGR}) but the method  is different and it can also
used to the fractional Brownian motion case and to the
multidimensional context.

\begin{prop}
Let $\xi$ be given by (\ref{xi}).  For every $x\in \mathbb{R}$,
it holds that
\begin{equation*}
\mathbf{E}\left|  \hat{\xi }(x) \right| ^{2} =T.
\end{equation*}
As a consequence $\xi \in H^{-r} (\mathbb{R}, \mathbb{R}) $  if and only if $r>\frac{1}{2}$.
\end{prop}
{\bf Proof: }We can write, by using the expression of the Skorohod
integral via chaos expansion
\begin{equation}\label{idelta1}
\xi (x)= \sum _{n\geq 0} I_{n+1}\left( \left( a_{n}^{x}(s)
1_{[0,s]}^{\otimes n}(\cdot ) \right) ^{\sim}\right)
\end{equation}
where $\left( a_{n}^{x}(s)
1_{[0,s]}^{\otimes n}(\cdot ) \right) ^{\sim}$
denotes the symmetrization in $n+1$ variables of the function
$(s,t_{1}, \ldots , t_{n})\to  a_{n}^{x}(s) 1_{[0,s]}^{\otimes
n}(t_{1}, \ldots , t_{n} ) $. It can be shown  that
\begin{equation*} \hat{\xi }(x) = \sum _{n\geq 0}I_{n+1}\left(
\left((a_{n}^{\hat{x}}(s) 1_{[0,s]}^{\otimes n}(\cdot ) \right) ^{\sim}
\right).
\end{equation*}
By using (\ref{imp}),
\begin{eqnarray*}
\hat{\xi }(x)&=&\sum _{n\geq 0}I_{n+1}\left( \left(
\frac{(-i)^{n}x^{n}}{n!} e^{-\frac{x^{2}s}{2}}1_{[0,s]}^{\otimes
n}(\cdot ) \right) ^{\sim} \right) =\sum _{n\geq
0}\frac{(-i)^{n}x^{n}}{n!}I_{n+1}\left( \left(
e^{-\frac{x^{2}s}{2}}1_{[0,s]}^{\otimes n}(\cdot ) \right) ^{\sim}\right)
\end{eqnarray*}
and  therefore, by the isometry of multiple integrals (\ref{iso}),
\begin{equation}
\label{rel2} \mathbf{E}\left|  \hat{\xi }(x) \right| ^{2}= \sum _{n\geq 0}
\frac{x^{2n}}{(n!)^{2}}(n+1)! \Vert \left(
e^{-\frac{x^{2}s}{2}}1_{[0,s]}^{\otimes n}(\cdot ) \right) ^{\sim} \Vert
^{2}_{L^{2}([0,T]^{n+1})}.
\end{equation}
Since
\begin{equation*}
\left( e^{-\frac{x^{2}s}{2}}1_{[0,s]}^{\otimes n}(\cdot )
\right) ^{\sim}(t_{1}, \ldots , t_{n+1}) = \frac{1}{n+1}\sum_{i=1} ^{n+1}
e^{-\frac{x^{2}t_{i}}{2}}1_{[0,t_{i}]}^{\otimes n}(t_{1}, \ldots ,
\hat{t_{i}}, \ldots , t_{n+1})
\end{equation*}
where $\hat{t_{i}}$ means that the variable $t_{i}$ is missing, we
obtain
\begin{eqnarray*}
\mathbf{E}\left|  \hat{\xi }(x) \right| ^{2}&=& \sum _{n\geq
0}\frac{x^{2n}}{(n+1)!} \Vert \sum_{i=0} ^{n+1}
e^{-\frac{x^{2}t_{i}}{2}}1_{[0,t_{i}]}^{\otimes n}(t_{1}, \ldots ,
\hat{t_{i}}, \ldots , t_{n+1}) \Vert ^{2} _{L^{2} ([0,T]^{n+1})} \\
&=&\sum _{n\geq 0}\frac{x^{2n}}{n!}\int_{0}^{T}ds \int _{[0,T]^{n}}
dt_{1}\ldots dt_{n} \left( e^{-\frac{x^{2}s}{2}}1_{[0,s]}^{\otimes
n}(t_{1}, \ldots , t_{n} )  \right) ^{2} \\
&=&\sum _{n\geq 0}\frac{x^{2n}}{n!}\int_{0}^{T}dse^{-x^{2}s}
s^{n}ds=\int_{0}^{T}e^{-x^{2}s} \sum_{n\geq
0}\frac{(x^{s}s)^{n}}{n!}ds= T.
\end{eqnarray*}
Finally,
\begin{equation*}
\mathbf{E}\left| \xi \right| _{H^{-r}(\mathbb{R}; \mathbb{R}) }^{2}= \int_{\mathbb{R}} (1+x^{2}
)^{-r}\mathbf{E}\left|  \hat{\xi }(x) \right| ^{2}dx=\int_{\mathbb{R}}
(1+x^{2} )^{-r}Tdx
\end{equation*}
and this is finite if and only if $r>\frac{1}{2}$. \qed

\vskip0.2cm

\begin{remark}
The above result gives actually a rigorous meaning of the formal calculation $E\left( \int_{0} ^{T} \delta (x-B_{s}) dB_{s} \right) ^{2} = E\int_{0}^{T} \vert \delta (x-B_{s}) \vert ^{2} ds = T$.
\end{remark}

\subsection{The multidimensional multiparameter Brownian currents}

Let us consider now the multidimensional situation. In this part we will actually  treat the regularity of the function
\begin{equation}
\label{xi2} \int _{[0,T] ^{N}} \delta (x-B_{s}  ) dB_{s}
\end{equation}
where $x\in \mathbb{R}^{d}$ and $B=\left( B^{1}, \ldots ,
B^{d}\right) $ is a $d$-dimensional Wiener sheet with parameter
$t\in [0,T] ^{N}$.  Here the term (\ref{xi2}) represents a vector
 given by
 \begin{eqnarray*}
\xi (x)&=&\int _{[0,T] ^{N}} \delta (x-B_{s}  ) dB_{s}= \left( \int
_{[0,T] ^{N}} \delta (x-B_{s}  ) dB^{1}_{s}, \ldots, \int _{[0,T]
^{N}} \delta (x-B_{s}  ) dB^{d}_{s} \right)\\
&:=& (\xi_{1}(x), \ldots , \xi_{d}(x))
 \end{eqnarray*}
 for every $x\in \mathbb{R}^{d}$.

The method considered above based on the chaos expansion of the delta Dirac function has  advantage that, in contrast with the approaches in  \cite{FGR} or \cite{FG}, it can be immediately extended to time parameter in $\mathbb{R}^{N}$. Moreover, we are able to compute explicitly the Sobolev norm of the vector  (\ref{xi2}) and to obtain an ``if and only if" result.
\begin{prop}
Let $(B_{t})_{t\in [0,T]^{N}}$ be a $d$ -dimensional Wiener process with multidimensional parameter $t\in [0,T]^{N}$. Then for every $w\in \Omega$ the mapping $x \to \int_{[0,T] ^{N}} \delta (x-B_{s}) dB_{s} $ belongs to the negative Sobolev space $H^{-r}(\mathbb{R}^{d}, \mathbb{R}^{d}) $ if and only if $r>\frac{d}{2}$.
\end{prop}
{\bf Proof: } We need to estimate
 $$\int_{\mathbb{R}^{d}} (1+\vert x\vert ^{2})^{-r} \mathbf{E}\left| \hat{\xi }(x) \right| ^{2}dx$$
with $ \left| \hat{\xi }(x) \right| ^{2}= \vert \hat{\xi} _{1}(x)\vert ^{2} + \ldots \vert \hat{\xi } _{d}(x)\vert ^{2}$
where $\vert  \hat{\xi }_{i}(x)  \vert$ denotes the complex modulus of  $\hat{\xi_{i}}(x)\in \mathbb{C}$.
\smallskip

We can formally write (but it can also written in a  rigorous manner by approximating the delta Dirac function by Gaussian kernels with variance $\eps \to 0$)
\begin{equation*}
\delta (x-B_{s}) =\prod _{j=1} ^{d} \delta (x_{j}-B_{s}^{j}) =
\delta _{k} (x-B_{s}) \delta (x_{k}-B_{s}^{k})
\end{equation*}
where we denoted
\begin{equation}
\label{deltak} \delta _{k} (x-B_{s}) = \prod _{j=1, j\not= k} ^{d}
\delta (x_{j}-B_{s}^{j}), \hskip0.5cm k=1, \ldots d.
\end{equation}
Let us compute the $k$th component of the vector $\xi (x)$. We will
use the chaotic expansion for the delta Dirac function with
multidimensional parameter $s\in [0,T] ^{N}$ (see (\ref{deltagen})
\begin{equation}
\label{deltaC} \delta (x_{j}-B_{s}^{j})= \sum _{n_{j}\geq 0}
a_{n_{j}}^{x_{j}} (s) I^{j}_{n_{j}}\left(
1_{[0,s]}^{\otimes n_{j}}(\cdot )\right)
\end{equation}
where $I^{j}_{n}$ denotes the Wiener-It\^o integral  of order $n$ with respect to the component $B^{j}$ and  for $s=(s_{1}, \ldots , s_{N})$
\begin{equation}
\label{anx2} a_{n_{j}}^{x_{j}} (s) =p_{s_{1}\ldots
s_{N}}(x_{j}) H_{n_{j}}\left( \frac{ x_{j}}{\sqrt{s_{1}\ldots
s_{N}}}\right) (s_{1} \ldots s_{N}) ^{-\frac{n_{j}}{2}}.
\end{equation}

It
holds, for every $k=1, \ldots , d$, by using formula (\ref{deltaB})
\begin{eqnarray*}
\xi_{k}(x)=\int _{[0,T] ^{N}} \delta (x-B_{s}  ) dB^{k}_{s}&=& \int
_{[0,T]^{N}}\delta _{k} (x-B_{s}) \delta
(x_{k}-B_{s}^{k})dB_{s}^{k}\\ &=&\int _{[0,T]^{N}}\delta _{k}
(x-B_{s})\sum _{n_{k}\geq 0} a_{n_{k}}^{x_{k}}(s)
I_{n_{k}}^{k} \left( 1_{[0,s]}^{\otimes n_{k}} (\cdot )\right)
dB_{s}^{k}.
\end{eqnarray*}
Since the components of the Brownian motion $B$ are independent, the
term $\delta _{k} (x-B_{s})$ is viewed as a deterministic function
when we integrate with respect to $B^{k}$. We obtain,
\begin{equation*}
\xi_{k}(x)=\int _{[0,T] ^{N}} \delta (x-B_{s}  ) dB^{k}_{s}=\sum _{n_{k}\geq 0}
I_{n_{k}+1} \left( \left(\delta _{k}(x-B_{s}) a
_{n_{k}}^{x_{k}}(s) 1_{[0,s]}^{\otimes n_{k}}(\cdot )\right) ^{\sim }\right).
\end{equation*}
Here $\left(\delta _{k}(x-B_{s}) a
_{n_{k}}^{x_{k}}(s) 1_{[0,s]}^{\otimes n_{k}}(\cdot )\right) ^{\sim }$ denotes the
symmetrization in $n_{k}+1$ variables of the function
$$(s,t_{1}, \ldots , t_{n_{k}}) =\delta _{k}(x-B_{s}) a
_{n_{k}}^{x_{k}}(s) 1_{[0,s]}^{\otimes n_{k}}(t_{1}, \ldots ,
t_{n_{k}}).$$

\vskip0.3cm

Let us denote by $\hat{ \xi } _{k}(x)$ the Fourier transform of $\xi
_{k}(x)= \int _{[0,T] ^{N}} \delta (x-B_{s}  ) dB^{k}_{s}$. As in
the previous section one can show that (that is, the Fourier
transform with respect to $x$ "goes inside" the  stochastic integral
)
\begin{equation}
\label{f1} \hat {\xi }_{k}(x)=\sum _{n_{k}\geq 0} I^{k}_{n_{k}+1} \left(
\left(\delta _{k}(\hat{x}-B_{s}) a
_{n_{k}}^{\hat{x_{k}}}(s) 1_{[0,s]}^{\otimes n_{k}}(\cdot )\right) ^{\sim} \right)
\end{equation}
where
$$\delta _{k}(\hat{x}-B_{s}) a
_{n_{k}}^{\hat{x_{k}}}(s) $$ is the Fourier transform of
$$\mathbb{R}^{d}\ni x =(x_{1}, \ldots , x_{d})\to\delta _{k}(x-B_{s}) a
_{n_{k}}^{x_{k}}(s)=\left( \prod _{j=1, j\not= k}^{d} \delta
(x_{j}-B_{s}^{j}) \right) a _{n_{k}}^{x_{k}}(s).$$ Now
clearly
\begin{equation*}
\delta _{k}(\hat{x}-B_{s}) =e^{-i\sum _{l=1, l\not= k} ^{d} x_{l}
B_{s}^{l}}
\end{equation*}
and by (\ref{imp})
\begin{equation*}
a _{n_{k}}^{\hat{x_{k}}}(s)=\frac{ (-i)^{n_{k}}x_{k}
^{n_{k}}}{n_{k}!} \mathbf{E} \left( e^{-ix_{k} B_{s}^{k}} \right)=\frac{ (-i)^{n_{k}}x_{k}
^{n_{k}}}{n_{k}!} e^{\frac{x_{k}^{2}\vert s\vert }{2}}
\end{equation*}
with $\vert s \vert = s_{1}\ldots s_{N}$ if $s=(s_{1}, \ldots , s_{N})$. Thus relation (\ref{f1}) becomes
\begin{equation}
\label{f2} \hat {\xi }_{k}(x)=\sum _{n_{k}\geq 0} I^{k}_{n_{k}+1} \left(
\left( e^{-i\sum _{l=1, l\not= k} ^{d} x_{l} B_{s}^{l}}\frac{
(-i)^{n_{k}}x_{k} ^{n_{k}}}{n_{k}!} \mathbf{E} \left( e^{-ix_{k} B_{s}^{k}}
\right) 1_{[0,s]}^{\otimes n_{k}}(\cdot )\right)^{\sim}\right)
\end{equation}
Taking the modulus in $\mathbb{C}$ of the above expression, we get
\begin{eqnarray*}
\vert \hat{\xi }_{k} (x)\vert ^{2}&=& \left| \sum _{n_{k}\geq 0} I^{k}_{n_{k}+1} \left( \left( \cos \left( \sum _{l=1, l\not= k} ^{d} x_{l} B_{s}^{l}  \right) \frac{
(-i)^{n_{k}}x_{k} ^{n_{k}}}{n_{k}!}e^{-\frac{x_{k}^{2}\vert s \vert }{2}}1_{[0,s]}^{\otimes n_{k}}(\cdot )\right)^{\sim }\right) \right|^{2}\\
&&+\left| \sum _{n_{k}\geq 0} I^{k}_{n_{k}+1} \left( \left( \sin \left( \sum _{l=1, l\not= k} ^{d} x_{l} B_{s}^{l}  \right) \frac{
(-i)^{n_{k}}x_{k} ^{n_{k}}}{n_{k}!}e^{-\frac{x_{k}^{2}\vert s \vert }{2}}1_{[0,s]}^{\otimes n_{k}}(\cdot )\right) ^{\sim} \right)\right|^{2}
\end{eqnarray*}
and using the isometry of multiple stochastic integrals (\ref{iso}) and the
independence of components we get as in the proof of Proposition 1,
(we use the notation $\vert s\vert =s_{1}\ldots s_{N}$)
\begin{eqnarray*}
\mathbf{E}\left| \hat{ \xi } _{k}(x) \right| ^{2}&=& \sum_{n_{k}\geq
0}\frac{x_{k} ^{2n_{k}}}{n_{k}!} \Vert  \cos \left( \sum _{l=1, l\not= k} ^{d} x_{l} B_{s}^{l}  \right) e^{-\frac{x_{k}^{2}\vert s \vert }{2}}1 _{[0,s]}
^{\otimes n_{k}} (\cdot ) \Vert ^{2} _{L^{2} ([0,T]^{N(n_{k}+1)} )
}\\
&&+ \sum_{n_{k}\geq
0}\frac{x_{k} ^{2n_{k}}}{n_{k}!} \Vert  \sin \left( \sum _{l=1, l\not= k} ^{d} x_{l} B_{s}^{l}  \right) e^{-\frac{x_{k}^{2}\vert s\vert }{2}}1 _{[0,s]}
^{\otimes n_{k}} (\cdot ) \Vert ^{2} _{L^{2} ([0,T]^{N(n_{k}+1)} )}\\
&=&\sum_{n_{k}\geq 0}\frac{x_{k} ^{2n_{k}}}{n_{k}!} \int
_{[0,T]^{N}}e^{ - x_{k} ^{2} \vert
s\vert } \vert s\vert ^{n_{k}} ds \\
&=&\int _{[0,T]^{N}} e^{ - x_{k} ^{2}\vert s\vert }\left( \sum_{n_{k}\geq 0}\frac{ (x_{k}^{2}
\vert s \vert )^{n_{k}}}{n_{k}!} \right)ds =\int _{[0,T]^{N}}  ds=T^{N}.
\end{eqnarray*}
This implies
\begin{eqnarray*}
\int_{\mathbb{R}^{d}}(1+\vert x\vert ^{2})^{-r}\mathbf{E}\left| \hat{ \xi}(x) \right| ^{2} dx&=&\int_{\mathbb{R}^{d}}(1+\vert x\vert ^{2})^{-r} \mathbf{E}\left(  \left|\hat{ \xi }
_{1}(x)
\right| ^{2}+ \ldots\left|  \hat{ \xi }_{d}(x) \right| ^{2}\right) dx\\
&=&dT^{N}\int_{\mathbb{R}^{d}}(1+\vert x\vert ^{2})^{-r}dx
\end{eqnarray*}
which is finite if and only if $2r>d$. \qed

\begin{remark}
It is interesting to observe that the dimension $N$ of the time parameter does not affect the regularity of currents. This is somehow unexpected because it is known that this dimension $N$ influences  the regularity of the local time of the process $B$ which can be formally written as $\int_{[0,T]^{N} }\delta (x-B_{s})ds$ (see \cite{Imk1}).
\end{remark}

\section{Regularity of  fractional  currents with respect to $x$}

\subsection{The one-dimensional fractional Brownian currents}

We will consider in this paragraph a fractional Brownian motion
$(B^{H}_{t}) _{t\in [0,T]}$ with  Hurst parameter $H\in
(\frac{1}{2}, 1)$. That is, $B^{H}$ is a centered Gaussian process
starting from zero with covariance function $$R^{H}(t,s):= \frac{1}{2}
\left( t^{2H}+s^{2H} -\vert t-s\vert ^{2H}\right), \hskip0.5cm s,t \in [0,T].$$ Let us denote
by ${\cal{H}}_{H}$ the canonical Hilbert space of the fractional
Brownian motion which is defined as the closure of the linear space
generated by the indicator functions $\{ 1_{[0,t] }, t\in [0,T]\}$
with respect to the scalar product
\begin{equation*}
\langle 1_{[0,t]}, 1_{[0,s]} \rangle _{{\cal{H}}_{H}} = R^{H}(t,s), \hskip0.5cm s,t\in [0,T].
\end{equation*}

One can construct multiple integrals with respect to $B^{H}$ (the
underlying space $L^{2}[0,T])$ is replaced by ${\cal{H}}_{H}$) and these integrals are those who appear in the formula (\ref{deltaB}).
We prove the following result.
\begin{prop} \label{p3} Let $B^{H}$ be a fractional Brownian motion with $H\in
(\frac{1}{2}, 1)$ and let $\xi $ be given by (\ref{xi}).  Then for every $w\in \Omega $, we have
\begin{equation*} \xi \in H^{-r} \left( \mathbb{R}; \mathbb{R}\right)
\end{equation*}
for every $r>\frac{1}{2H}-\frac{1}{2}$.
\end{prop}
{\bf Proof: }
  Using
the above computations, we will get that the Fourier transform of
$g(x)=\delta(x-B^{H}_{s} ) $ is equal, on one hand, to
\begin{equation*}
\hat{g} (x)= \sum_{n\geq 0}a_{n}^{\hat{x}} (s^{2H})I_{n}^{B^{H}}
\left( 1_{[0,t]}^{\otimes n}(\cdot )\right)
\end{equation*}
(here $I_{n}^{B^{H}}$ denotes the multiple stochastic integral with
respect to the fBm $B^{H}$ and the function $a_{n}^{x}$ is defined
by (\ref{an})), and on the other hand,

\begin{equation*}
\hat{g} (x)= \mathbf{E}\left(  e^{-ixB^{H}_{s} }\right)\sum _{n\geq
0}\frac{(-ix)^{n} }{n!} I_{n}^{B} \left( 1_{[0,t]}^{\otimes n}(\cdot
)\right)
\end{equation*}
and so
\begin{equation}
\label{imp2} a_{n}^{\hat{x}} (s^{2H})= e^{-\frac{x^{2} s^{2H}
}{2}}\frac{(-ix)^{n} }{n!}.
\end{equation}
Moreover
\begin{equation*}
\hat{\xi }(x)= \sum _{n\geq
0}\frac{(-i)^{n}x^{n}}{n!}I^{B^{H}}_{n+1}\left( \left(
e^{-\frac{x^{2}s^{2H}}{2}} 1_{[0,s]}^{\otimes n}(\cdot ) \right) ^{\sim}\right).
\end{equation*}
We then obtain
\begin{eqnarray*}
\mathbf{E}\Vert \xi \Vert ^{2} _{H^{-r}(\mathbb{R}; \mathbb{R})}&=&\mathbf{E} \int_{\mathbb{R}}
\frac{1}{(1+x^{2})^{r} }\left|\hat{\xi }(x) \right| ^{2}dx\\
&=&\int_{\mathbb{R}} \frac{1}{(1+x^{2})^{r} }\sum _{n\geq 0}
\frac{x^{2n}}{(n!)^{2}}(n+1)! \Vert \left(
e^{-\frac{x^{2}s^{2H}}{2}}1_{[0,s]}^{\otimes n}(\cdot ) \right) ^{\sim} \Vert
^{2}_{{\cal{H}}^{\otimes n+1}}dx\\
&=&\int_{\mathbb{R}} \frac{1}{(1+x^{2})^{r} }\sum _{n\geq 0}
\frac{x^{2n}}{(n+1)!}dx\\
&&\times \langle \sum
_{i=1}^{n+1}e^{-\frac{x^{2}u_{i}^{2H}}{2}}1_{[0,u_{i}]}^{\otimes
n}(u_{1}, \ldots , \hat{u}_{i}, \ldots , u_{n+1}), \sum
_{j=1}^{n+1}e^{-\frac{x^{2}v_{j}^{2H}}{2}}1_{[0,v_{j}]}^{\otimes
n}(v_{1}, \ldots , \hat{v}_{j}, \ldots , v_{n+1})\rangle
_{{\cal{H}}_{H}}.
\end{eqnarray*}
Using the fact that for regular enough functions  $f$ and $g$ in
${\cal{H}}_{H}$ their product scalar is given by (see e.g. \cite{N},
Chapter 5) $ \langle f,g \rangle  _{{\cal{H}}_{H}}= H(2H-1) \int_{0}^{T}
\int_{0}^{T} f(x)g(y) \vert x-y\b\vert ^{2H-2} dxdy $ and thus for
regular enough function $f,g\in {\cal{H}}^{\otimes n}$
\begin{equation*}
\langle f,g \rangle  _{{\cal{H}}_{H}^{\otimes n}}=(H(2H-1))^{n}
\int_{[0,T]^{n}}\int_{[0,T]^{n}}f(u_{1}, \ldots , u_{n}) g(v_{1},
\ldots , v_{n})\prod_{i=1}^{n} \vert u_{i} -v_{i}\vert ^{2H-2}du_{1}..du_{n}dv_{1}..dv_{n}
\end{equation*}
we find (by $du_{i}du_{j}$ we mean below $du_{1}..du_{n+1}dv_{1}..dv_{n+1}$)

\begin{eqnarray*}
&&\mathbf{E}\Vert \xi \Vert _{H^{-r}(\mathbb{R}; \mathbb{R})}^{2}\\
&=&\int_{\mathbb{R}}
\frac{1}{(1+x^{2})^{r} }\sum _{n\geq 0}(H(2H-1))^{n+1}
\frac{x^{2n}}{(n+1)!}dx  \sum _{i,j=1}^{n+1} \int _{[0,T]^{n+1}}
\int _{[0,T]^{n+1}} \prod_{l=1}^{n+1} \vert u_{l}-v_{l} \vert
^{2H-2} du_{i}dv_{j}\\
&&\times
e^{-\frac{x^{2}u_{i}^{2H}}{2}}e^{-\frac{x^{2}v_{j}^{2H}}{2}}1_{[0,u_{i}]}^{\otimes
n}(u_{1}, \ldots , \hat{u}_{i}, \ldots ,
u_{n+1})1_{[0,v_{j}]}^{\otimes n}(v_{1}, \ldots , \hat{v}_{j},
\ldots , v_{n+1})\\
&=& \int_{\mathbb{R}} \frac{1}{(1+x^{2})^{r} }\sum _{n\geq
0}(H(2H-1))^{n+1} \frac{x^{2n}}{(n+1)!}dx \\
&&\times
 \sum _{i=1}^{n+1} \int _{[0,T]^{n+1}} \int _{[0,T]^{n+1}}
e^{-\frac{x^{2}u_{i}^{2H}}{2}}e^{-\frac{x^{2}v_{i}^{2H}}{2}}\prod_{l=1}^{n+1}
\vert u_{l}-v_{l} \vert ^{2H-2} du_{i}dv_{j}\\
&+& \int_{\mathbb{R}} \frac{1}{(1+x^{2})^{r} }\sum _{n\geq
1}(2H(2H-1))^{n+1 }\frac{x^{2n}}{(n+1)!}dx\\
&&\times \sum _{i,j=1; i\not=j}^{n+1} \int _{[0,T]^{n+1}} \int
_{[0,T]^{n+1}} \prod_{l=1}^{n+1} \vert u_{l}-v_{l} \vert
^{2H-2} du_{i}dv_{j}\\
&&\times
e^{-\frac{x^{2}u_{i}^{2H}}{2}}e^{-\frac{x^{2}v_{j}^{2H}}{2}}1_{[0,u_{i}]}^{\otimes
n}(u_{1}, \ldots , \hat{u}_{i}, \ldots ,
u_{n+1})1_{[0,v_{j}]}^{\otimes n}(v_{1}, \ldots , \hat{v}_{j},
\ldots , v_{n+1})\\
&:=& A+B.
\end{eqnarray*}
Let us compute first the term $A$. Using the symmetry of the
integrand and the fact that
\begin{equation}\label{II}
2H(2H-1)\int_{0}^{t}\int_{0}^{s} \vert u-v\vert ^{2H-2} dudv=
R(t,s), \hskip0.5cm s,t\in [0,T]
\end{equation}
we can write
\begin{eqnarray*}
A&=&\int_{\mathbb{R}} \frac{1}{(1+x^{2})^{r} }\sum _{n\geq
0}(2H(2H-1))^{n+1} \frac{x^{2n}}{n!}dx \\
&&\times \sum _{i=1}^{n+1} \int _{[0,T]^{n+1}} \int _{[0,T]^{n+1}}
\prod_{l=1}^{n+1} \vert u_{l}-v_{l} \vert
^{2H-2} du_{i}dv_{j}\\
&&\times
e^{-\frac{x^{2}u_{1}^{2H}}{2}}e^{-\frac{x^{2}v_{1}^{2H}}{2}}1_{[0,u_{1}]}^{\otimes
n}(u_{2}, \ldots ,, u_{n+1})1_{[0,v_{1}]}^{\otimes n}(v_{2},
\ldots v_{n+1})
\end{eqnarray*}
and by integrating with respect to $du_{2} \ldots d u_{n+1} $.
\begin{eqnarray*}
A&=&\int_{\mathbb{R}}dx \frac{1}{(1+x^{2})^{r} }\sum _{n\geq
0}2H(2H-1) \frac{x^{2n}}{n!}dx \\
&&\times \int_{0}^{T} \int
_{0}^{T}e^{-\frac{x^{2}u^{2H}}{2}}e^{-\frac{x^{2}v^{2H}}{2}}\vert
u-v\vert ^{2H-2}R(u,v)^{n} dudv \\
&=& H(2H-1)\int_{\mathbb{R}}dx \frac{1}{(1+x^{2})^{r} } \int
_{0}^{T} \int
_{0}^{T}e^{-\frac{x^{2}u^{2H}}{2}}e^{-\frac{x^{2}v^{2H}}{2}}\vert
u-v\vert ^{2H-2} \left( \sum _{n\geq 0}
\frac{x^{2n}R(u,v)^{n}}{n!}\right)dudv\\
&=& H(2H-1)\int_{\mathbb{R}} dx\frac{1}{(1+x^{2})^{r} } \int
_{0}^{T} \int
_{0}^{T}e^{-\frac{x^{2}u^{2H}}{2}}e^{-\frac{x^{2}v^{2H}}{2}}\vert
u-v\vert ^{2H-2} e^{x^{2} R(u,v)} dudv
\\
&=&H(2H-1)\int_{\mathbb{R}} \frac{1}{(1+x^{2})^{r} } \int _{0}^{T}
\int _{0}^{T}e^{-\frac{x^{2}\vert u-v\vert ^{2H}}{2}}\vert
u-v\vert ^{2H-2} dudv.
\end{eqnarray*}
Now, by classical Fubini,
\begin{equation*}
A=H(2H-1) \int _{0}^{T} \int _{0}^{T}\vert u-v\vert ^{2H-2}
dudv\int_{\mathbb{R}} \frac{1}{(1+x^{2})^{r} }e^{-\frac{x^{2}\vert
u-v\vert ^{2H}}{2}}dx
\end{equation*}
and by using the change of variable $ x\vert u-v\vert ^{H}=y $ in
the integral with respect to ,
\begin{equation*}
A=H(2H-1) \int _{0}^{T} \int _{0}^{T}\vert u-v\vert
^{2H-2+2Hr-H}\int_{\mathbb{R}} e^{-\frac{y^{2}}{2}} \left( \vert
u-v\vert ^{2H} + y^{2} \right) ^{-r} dy
\end{equation*}
and this is finite when $2H-2+2Hr-H >-1$ which gives
$r>\frac{1}{2H}-\frac{1}{2}$.

\vskip0.5cm

The term denoted by $B$ can be treated as follows.
\begin{eqnarray*}
B&=&\int_{\mathbb{R}}\ \frac{1}{(1+x^{2})^{r} }\sum _{n\geq 1}
\frac{x^{2n}}{(n+1)!} n(n+1) (2H(2H-1)) ^{n+1}\int _{[0,T]^{n+1}}
\int _{[0,T]^{n+1}} du_{i}dv_{j} \\
&&\times e^{-\frac{x^{2} u_{1}^{2H}}{2}}e^{-\frac{x^{2}
v_{2}^{2H}}{2}}1_{[0, u_{1} ]}(u_{2},\ldots , u_{n+1})^{\otimes n}1_{[0, v_{2}
]}(v_{1}, v_{3},\ldots , v_{n+1})^{\otimes n}\prod_{l=1}^{n+1} \vert
u_{l}-v_{l} \vert ^{2H-2}\\
&=&(2H(2H-1))^{2}\int_{\mathbb{R}}\ \frac{1}{(1+x^{2})^{r} }\sum
_{n\geq 1} \frac{x^{2n}}{(n-1)!}\int _{0}^{T} \int _{0}^{T}\int
_{0}^{T} \int _{0}^{T}du_{1}du_{2}dv_{1}dv_{2}\\
&&\times e^{-\frac{x^{2} u_{1}^{2H}}{2}}e^{-\frac{x^{2}
v_{2}^{2H}}{2}}1_{[0, u_{1} ]}(u_{2})1_{[0,v_{2}]}(v_{1}) R(u_{1},
v_{2}) ^{n-1}\vert u_{1}-v_{1} \vert ^{2H-2} \vert u_{2}-v_{2}
\vert ^{2H-2}.
\end{eqnarray*}
Since
\begin{equation*}
\sum _{n\geq 1} \frac{x^{2n}}{(n-1)!}R(u_{1}, v_{2}) ^{n-1}= x^{2}
\sum _{n\geq 0} \frac{x^{2n} R(u_{1}, v_{2})^{n}}{n!} =x^{2}
e^{x^{2} R(u_{1}, v_{2} )}
\end{equation*}
we obtain
\begin{eqnarray*}
B&=&(2H(2H-1))^{2}\int_{\mathbb{R}}\ \frac{1}{(1+x^{2})^{r}
}x^{2}\int _{0}^{T} \int _{0}^{T}\int
_{0}^{T} \int _{0}^{T}du_{1}du_{2}dv_{1}dv_{2}\\
&&\times e^{-\frac{x^{2}\vert u_{1}-v_{2}\vert ^{2H}}{2}}\vert
u_{1}-v_{1} \vert ^{2H-2} \vert u_{2}-v_{2} \vert ^{2H-2}1_{[0,
u_{1} ]}(u_{2})1_{[0,v_{2}]}(v_{1}).
\end{eqnarray*}
We calculate first the integral $du_{2}$ and $dv_{1}$ and then in
$dx$ we use the change of variables $x\vert u_{1}-v_{2}\vert
^{H}=y$ and we get
\begin{eqnarray*}
B&=& H^{2}\int_{\mathbb{R}} \frac{1}{(1+x^{2})^{r} }x^{2}\int
_{0}^{T} \int _{0}^{T} du_{1}dv_{2}e^{-\frac{x^{2}\vert
u_{1}-v_{2}\vert
^{2H}}{2}}\\
&&\times \left( \vert u_{1} - v_{2}\vert ^{2H-1}  -
u_{1}^{2H-1}\right)\left( \vert u_{1} - v_{2}\vert ^{2H-1}  -
v_{2}^{2H-1}\right)
\\
&=& \int _{0}^{T} du_{1}dv_{2}\left( \vert u_{1} - v_{2}\vert
^{2H-1}  - u_{1}^{2H-1}\right)\left( \vert u_{1} - v_{2}\vert
^{2H-1}  - v_{2}^{2H-1}\right)\vert u_{1}-v_{2} \vert
^{2Hr-2H-H}\\
&&\times
\int_{\mathbb{R}}y^{2} e^{-\frac{y^{2}}{2}} \left( \vert
u_{1}-v_{2} \vert ^{2H} +y^{2} \right) dy.
\end{eqnarray*}
which is finite if $4H-2 +2Hr-3H > -1$ and this implies again
$r>\frac{1}{2H}-\frac{1}{2}$. \qed

 \vskip0.5cm

 \subsection{The multidimensional fractional currents}

We will also consider the multidimensional case of the fractional  Brownian motion $B^{H}$ in $\mathbb{R}^{d}$. It is defined as a random vector $B^{H}=(B^{H_{1}}, \ldots ,B^{H_{d}})$ where $B^{H_{i}}$ are independent one-dimensional fractional Brownian sheets. We will assume that  $H_{i} \in (\frac{1}{2}, 1)$ for every $i=1, \ldots , d$.  In this case
\begin{equation*}
\xi (x)= \left( \int_{[0,T] ^{N}} \delta (x-B^{H}_{s}) dB^{H_{1}}, \ldots ,\int_{[0,T] ^{N}} \delta (x-B^{H}_{s}) dB^{H_{d}}\right).
\end{equation*}
where the above integral is a divergence (Skorohod) integral with respect to the fractional Brownian motion.  We mention that the canonical Hilbert space ${\cal{H}}_{H_{k}}$ of the fractional Brownian sheet $B^{H_{k}}$ is now the closure of the linear space of the indicator functions with respect to the inner product
\begin{equation*}
\langle 1_{[0,t] }, 1_{[0,s]}\rangle _{{\cal{H}}_{H_{k}}}=\mathbf{E}(B^{H_{k}}_{t}B^{H_{k}}_{s}):= R^{H_{k}}(t,s)= \prod _{i=1}^{N} R^{H_{k}}(t_{i}, s_{i})
\end{equation*} if $t=(t_{1}, \ldots , t_{N})$ and $s=(s_{1}, \ldots , s_{N})$. In this case we have

 \begin{prop}
 Let $B^{H}$ be as above. Then for every $\omega$ the fractional Brownian current $\xi $ belong to the Sobolev space $H^{-r}(\mathbb{R}^{d}; \mathbb{R}^{d})$ for every $r>\max _{k=1,..,d}\left( \frac{d}{2}-1+\frac{1}{2H_{k}}\right)$.
 \end{prop}
{\bf Proof: }We write \begin{eqnarray*}
\xi_{k}(x)=\int _{[0,T] ^{N}} \delta (x-B^{H}_{s}  ) dB^{H_{k}}_{s}&=& \int
_{[0,T]^{N}}\delta _{k} (x-B^{H}_{s}) \delta
(x_{k}-B^{H_{k}}_{s})dB_{s}^{H_{k}}\\ &=&\int _{[0,T]^{N}}\delta _{k}
(x-B^{H}_{s})\sum _{n_{k}\geq 0} a_{n_{k}}^{x_{k}}(s)
I_{n_{k}}^{H_{k}} \left( 1_{[0,s]}^{\otimes n_{k}} (\cdot )\right)
dB_{s}^{H_{k}}.
\end{eqnarray*}
and\begin{equation*}
 \hat {\xi }_{k}(x)=\sum _{n_{k}\geq 0} I^{H_{k}}_{n_{k}+1} \left( \left(
 e^{-i\sum _{l=1, l\not= k} ^{d} x_{l} B_{s}^{H_{l}}}\frac{
(-i)^{n_{k}}x_{k} ^{n_{k}}}{n_{k}!} e^{-\frac{x_{2}\vert s\vert ^{2H_{k}}}{2}}
 1_{[0,s]}^{\otimes n_{k}}(\cdot )\right) ^{\sim }\right)
\end{equation*}
Here again the integral $dB^{H_{k}}$ denotes the Skorohod integral with respect to the fractional Brownian sheet $B^{H_{k}}$ and $I^{H_{k}}$ denotes the multiple integral with respect to $B^{H_{k}}$. Then
\begin{eqnarray*}
\left| \hat {\xi }_{k}(x)\right| ^{2}&=& \left|\sum _{n_{k}\geq 0} I^{H_{k}}_{n_{k}+1} \left(\left( \cos \left( \sum _{l=1, l\not= k} ^{d} x_{l} B_{s}^{H_{l}}\right) \frac{
x_{k} ^{n_{k}}}{n_{k}!} e^{-\frac{x^{2}\vert s\vert ^{2H_{k}}}{2}}
 1_{[0,s]}^{\otimes n_{k}}(\cdot )\right) ^{\sim }\right)\right| ^{2}\\
 &&+\left|\sum _{n_{k}\geq 0} I^{H_{k}}_{n_{k}+1} \left(\left( \sin \left( \sum _{l=1, l\not= k} ^{d} x_{l} B_{s}^{H_{l}}\right) \frac{
x_{k} ^{n_{k}}}{n_{k}!} e^{-\frac{x^{2}\vert s\vert ^{2H_{k}}}{2}}
 1_{[0,s]}^{\otimes n_{k}}(\cdot )\right) ^{\sim }\right)\right| ^{2}
 \end{eqnarray*}
 and
 \begin{eqnarray*}
 \mathbf{E}\left| \hat {\xi }_{k}(x)\right| ^{2}&=&\sum_{n_{k}\geq 0} (n_{k}+1)!\frac{
x_{k} ^{2n_{k}}}{(n_{k}!)^{2}} \left| \left| \left( \cos \left( \sum _{l=1, l\not= k} ^{d} x_{l} B_{s}^{H_{l}}\right)  e^{-\frac{x^{2}\vert s\vert ^{2H_{k}}}{2}}
 1_{[0,s]}^{\otimes n_{k}}(\cdot )\right) ^{\sim }\right| \right| ^{2} _{({\cal{H}}_{H_{k}})^{\otimes (n_{k}+1)}}\\
 &&+\sum_{n_{k}\geq 0} (n_{k}+1)! \frac{
x_{k} ^{2n_{k}}}{(n_{k}!)^{2}} \left| \left| \left( \sin \left( \sum _{l=1, l\not= k} ^{d} x_{l} B_{s}^{H_{l}}\right) e^{-\frac{x^{2}\vert s\vert ^{2H_{k}}}{2}}
 1_{[0,s]}^{\otimes n_{k}}(\cdot )\right) ^{\sim }\right| \right| ^{2} _{({\cal{H}}_{H_{k}})^{\otimes (n_{k}+1)}}
 \end{eqnarray*}
 Note that, if $f,g$ are two regular functions of $N$ variables
 \begin{eqnarray*}
 \langle f,g\rangle _{{\cal{H}}_{H_{k}}}&=&(H_{k}(2H_{k}-1)) ^{N} \int_{[0,T]^{N}}\int_{[0,T]^{N}}f(u)g(u) \prod_{i=1}^{N}\vert u_{i}-v_{i}\vert ^{2H_{k}-2}dudv\\
 &:=&(H_{k}(2H_{k}-1)) ^{N} \int_{[0,T]^{N}}\int_{[0,T]^{N}}f(u)g(v)\Vert u-v\Vert ^{2H_{k}-2}dudv
 \end{eqnarray*}
if $u=(u_{1}, ..,u_{N})$ and $v=(v_{1},..,v_{N})$ and if $F,G$ are regular functions of $N(n_{k}+1)$ variables,
\begin{eqnarray*}
&&\langle F,G\rangle  _{\left({\cal{H}}_{H_{k}}\right) ^{\otimes (n_{k}+1)}}=\prod_{k=1}^{n_{k}+1}(H_{k}(2H_{k}-1)) ^{N} \\
&&\int_{[0,T]^{N(n_{k}+1)}}\int_{[0,T]^{N(n_{k}+1)}}F(u^{1}, .., u^{n_{k}+1})G(v^{1}, ..,v^{n_{k}+1})\prod_{i=1}^{n_{k}+1}\Vert u^{i}-v^{i}\Vert ^{2H_{k}-2}du^{1}..du^{n_{k}+1}dv^{1}..dv^{n_{k}+1}
\end{eqnarray*}
and thus, by symmetrizing the above function and taking the scalar product in $\left({\cal{H}}_{H_{k}}\right) ^{\otimes (n_{k}+1)}$
\begin{eqnarray*}
&&\mathbf{E}\left| \hat {\xi }_{k}(x)\right| ^{2}\\
&=&\sum_{n_{k}\geq 0}\frac{x_{k}^{2n_{k}}}{(n_{k}+1)!}(H_{k}(2H_{K}-1))^{N(n_{k}+1)} \sum_{i,j=1}^{n_{k}+1} \int_{[0,T]^{N(n_{k}+1)}}du\int_{[0,T]^{N(n_{k}+1)}}dv\prod_{l=1}^{n_{k}+1} \Vert u^{l}-v^{l}\Vert ^{2H_{k}-2}\\
&&\times  \cos \left( \sum _{l=1, l\not= k} ^{d} x_{l} B_{u^{i}}^{H_{l}}\right) \cos \left( \sum _{l=1, l\not= k} ^{d} x_{l} B_{v^{j}}^{H_{l}}\right)
e^{-\frac{x_{k}^{2}\vert u^{i}\vert ^{2H_{k}}}{2}}e^{-\frac{x_{k}^{2}\vert v^{j}\vert ^{2H_{k}}}{2}}\\
&&\times 1_{[0,u^{i}]}^{\otimes
n_{k}}(u^{1}, \ldots , \hat{u}^{i}, \ldots ,
u^{n_{k}+1})1_{[0,v^{j}]}^{\otimes n_{k}}(v^{1}, \ldots , \hat{v}^{j},
\ldots , v^{n_{k}+1})\\
&+&\sum_{n_{k}\geq 0}\frac{x_{k}^{2n_{k}}}{(n_{k}+1)!}(H_{k}(2H_{K}-1))^{N(n_{k}+1)} \sum_{i,j=1}^{n_{k}+1} \int_{[0,T]^{N(n_{k}+1)}}du\int_{[0,T]^{N(n_{k}+1)}}dv\prod_{l=1}^{n_{k}+1} \Vert u^{l}-v^{l}\Vert ^{2H_{k}-2}\\
&&\times  \sin \left( \sum _{l=1, l\not= k} ^{d} x_{l} B_{u^{i}}^{H_{l}}\right) \sin \left( \sum _{l=1, l\not= k} ^{d} x_{l} B_{v^{j}}^{H_{l}}\right)
e^{-\frac{x_{k}^{2}\vert u^{i}\vert ^{2H_{k}}}{2}}e^{-\frac{x_{k}^{2}\vert v^{j}\vert ^{2H_{k}}}{2}}\\
&&\times 1_{[0,u^{i}]}^{\otimes
n_{k}}(u^{1}, \ldots , \hat{u}^{i}, \ldots ,
u^{n_{k}+1})1_{[0,v^{j}]}^{\otimes n_{k}}(v^{1}, \ldots , \hat{v}^{j},
\ldots , v^{n_{k}+1})
\end{eqnarray*}
Here $\vert u^{i}\vert ^{2H_{k}}=(u^{i}_{1}...u^{i}_{N})^{2H_{k}} $ and $1_{[0, u^{i}]}= 1_{[0, u^{i}_{1}]}..1_{[0, u^{i}_{N}]}$ if $u^{i}=(u^{i}_{1}, .., u^{i}_{N})$. The next step is to majorize $\vert (\cos(u)\cos (v)\vert $ by $\frac{1}{2} (\cos^{2} (u)+ \cos^{2}(v))$ and similarly for the sinus.
\begin{eqnarray*}
&&\mathbf{E}\left| \hat {\xi }_{k}(x)\right| ^{2}\\
&\leq &cst. \sum_{n_{k}\geq 0}\frac{x_{k}^{2n_{k}}}{(n_{k}+1)!}(H_{k}(2H_{K}-1))^{N(n_{k}+1)} \sum_{i=1}^{n_{k}+1} \int_{[0,T]^{N(n_{k}+1)}}du\int_{[0,T]^{N(n_{k}+1)}}dv\prod_{l=1}^{n_{k}+1} \Vert u^{l}-v^{l}\Vert ^{2H_{k}-2}\\
&&e^{-\frac{x^{2}\vert u^{i}\vert ^{2H_{k}}}{2}}e^{-\frac{x^{2}\vert v^{i}\vert ^{2H_{k}}}{2}} 1_{[0,u^{i}]}^{\otimes
n_{k}}(u^{1}, \ldots , \hat{u}^{i}, \ldots ,
u^{n_{k}+1})1_{[0,v^{j}]}^{\otimes n_{k}}(v^{1}, \ldots , \hat{v}^{j},
\ldots , v^{n_{k}+1})\\
&&+cst. \sum_{n_{k}\geq 0}\frac{x_{k}^{2n_{k}}}{(n_{k}+1)!}(H_{k}(2H_{K}-1))^{N(n_{k}+1)} \sum_{i,j=1;i\not= j}^{n_{k}+1} \int_{[0,T]^{N(n_{k}+1)}}du\int_{[0,T]^{N(n_{k}+1)}}dv\prod_{l=1}^{n_{k}+1} \Vert u^{l}-v^{l}\Vert ^{2H_{k}-2}\\
&&e^{-\frac{x_{k}^{2}\vert u^{i}\vert ^{2H_{k}}}{2}}e^{-\frac{x_{k}^{2}\vert v^{j}\vert ^{2H_{k}}}{2}} 1_{[0,u^{i}]}^{\otimes
n_{k}}(u^{1}, \ldots , \hat{u}^{i}, \ldots ,
u^{n_{k}+1})1_{[0,v^{j}]}^{\otimes n_{k}}(v^{1}, \ldots , \hat{v}^{j},
\ldots , v^{n_{k}+1}):=C_{k}+D_{k}
\end{eqnarray*}
and the two terms above can be treated as the terms $A$ and $B$ from the one-dimensional case. For example the term denotes by $C_{k}$ (we illustrate only this term because $D_{k}$ is similar to the term $B$ in the one-dimensional case) gives
\begin{eqnarray*}
C_{k}&=&\sum_{n_{k}\geq 0}\frac{x_{k}^{2n_{k}}}{n_{k}!}(H_{k}(2H_{K}-1))^{N(n_{k}+1)}\int_{[0,T]^{N(n_{k}+1)}}\int_{[0,T]^{N(n_{k}+1)}}\prod_{l=1}^{n_{k}+1} \Vert u^{l}-v^{l}\Vert ^{2H_{k}-2}\\
&&e^{-\frac{x_{k}^{2}\vert u^{1}\vert ^{2H_{k}}}{2}}e^{-\frac{x_{k}^{2}\vert v^{1}\vert ^{2H_{k}}}{2}} 1_{[0,u^{1}]}^{\otimes
n_{k}}(u^{2}, \ldots ,  \ldots ,
u^{n_{k}+1})1_{[0,v^{1}]}^{\otimes n_{k}}(v^{2}, \ldots ,
\ldots , v^{n_{k}+1})\\
&=&\sum_{n_{k}\geq 0}\frac{x_{k}^{2n_{k}}}{n_{k}!}(H_{k}(2H_{k}-1))^{N}\int_{[0,T]^{N}}\int_{[0,T]^{N}}e^{-\frac{x_{k}^{2}\vert u^{1}\vert ^{2H_{k}}}{2}} e^{-\frac{x_{k}^{2}\vert v^{1}\vert ^{2H_{k}}}{2}}\Vert u^{1}-v^{1}\Vert ^{2H_{k}-2}R^{H_{k}}(u^{1}, v^{1}) ^{n_{k}}d^{1}dv^{1}\\
&=&(H_{k}(2H_{k}-1))^{N}\int_{[0,T]^{N}}\int_{[0,T]^{N}}e^{x_{k}^{2}R(u^{1}, v^{1})}e^{-\frac{x_{k}^{2}\vert u^{1}\vert ^{2H_{k}}}{2}} e^{-\frac{x_{k}^{2}\vert v^{1}\vert ^{2H_{k}}}{2}}\Vert u^{1}-v^{1}\Vert ^{2H_{k}-2}d^{1}dv^{1}
\end{eqnarray*}
with $R^{H_{k}}(u^{1}, v^{1})=R^{H_{k}}(u^{1}_{1},v^{1}_{1}) ..R^{H_{k}}(u^{1}_{N},v^{1}_{N})$.
But $R^{H_{k}}(u,v)-u^{2H_{k}}-v^{2H_{k}}=\mathbf{E}(B_{u}^{H_{k}}B_{v}^{H_{k}})-\mathbf{E}(B^{H_{k}}_{u})^{2}-\mathbf{E}(B^{H_{k}}_{v})^{2}=-\mathbf{E}(B^{H_{k}}_{u}-B^{H_{k}}_{v})^{2}$ for any $u,v\in \mathbb{R}^{n}$ and this implies
$$\int_{\mathbb{R}^{d}}dx(1+\vert x\vert ^{2})^{-r}C_{k}=cst. \int_{\mathbb{R}^{d}}dx(1+\vert x\vert ^{2})^{-r}\int_{[0,T]^{N}}\int_{[0,T]^{N}}e^{-\frac{x_{k}^{2}\mathbf{E}(B^{H_{k}}_{u}-B^{H_{k}}_{v})^{2}}{2}}\Vert u-v\Vert ^{2H_{k}-2}dudv.$$
The case $N=1$ can be easily handled. Indeed, in this case
\begin{eqnarray*}
\int_{\mathbb{R}^{d}}dx(1+\vert x\vert ^{2})^{-r}C_{k}&=&
\int_{[0,T]}\int_{[0,T]} \vert u-v\vert ^{2H_{k}-2}\int_{\mathbb{R}^{d}} (1+\vert x\vert )^{-r}e^{-\frac{x_{k}^{2}\vert u-v\vert ^{2H_{k}}}{2}}dx\\
&=&\int_{[0,T]}\int_{[0,T]} \vert u-v\vert ^{2H_{k}-2+2H_{k}r-dH_{k}}\int_{\mathbb{R}^{d}}e^{-\frac{y_{k}^{2}}{2}}(\vert u-v\vert ^{2H_{k}}+y^{2})^{-r}dy
\end{eqnarray*}
which is finite if $r>\frac{1}{2H_{k}}+\frac{d}{2}-1$. When $N\geq 2$ we will have
 \begin{eqnarray*}
&&\int_{\mathbb{R}^{d}}dx(1+\vert x\vert ^{2})^{-r}C_{k}\\
&=&\int_{[0,T]^{N}}\int_{[0,T]^{N}}\Vert u-v\Vert ^{2H_{k}-2}\left( \mathbf{E}(B^{H_{k}}_{u}-B^{H_{k}}_{v})^{2}\right) ^{\frac{d}{2}+r}\int_{\mathbb{R}^{d}}e^{-y_{k}^{2}/2} \left[ \mathbf{E}(B^{H_{k}}_{u}-B^{H_{k}}_{v})^{2}+ y^{2}\right] ^{-r}dx
\end{eqnarray*}
and this is finite if $r>\frac{1}{2H_{k}}+\frac{d}{2}-1$.
\begin{remark}
In \cite{FGR} the authors obtained the same regularity in the case of the pathwise integral with respect to the fBm. They  considered only the case $H_{k}=H$ for every $k=1,.,d$ and $N=1$. On the other hand, \cite{FGR} the Hurst parameter is allowed to be lesser than one half, it is assumed to be in $(\frac{1}{4}, 1)$.
\end{remark}

\section{Regularity of stochastic currents  with respect to $\omega$}

\subsection{The Brownian case}

We study now, for fixed $x\in \mathbb{R}$, the regularity of the
functional $\xi$ in the sense of Watanabe. We would like to see if, as for the delta Dirac function, the stochastic integral $\xi (x)$ keeps the same regularity in the Watanabe spaces (as a function of $\omega$) and in the Sobolev spaces as a function of $x$. We consider $B$ a one -dimensional Wiener process and denote
$$\xi (x)= \int _{a}^{T} \delta (x-B_{s} ) dB_{s}$$ with $a>0$.
Recall that by (\ref{deltaB})
\begin{equation*}\label{idelta}
\xi (x)= \sum _{n\geq 0} I_{n+1}\left( \left(a_{n}^{x}(s)
1_{[0,s]}^{\otimes n}(\cdot ) 1_{[a,T]}(s)\right) ^{\sim}\right)
\end{equation*}
(we consider the integral from $a>0$ instead of zero to avoid a singularity) where
\begin{equation*}
\left(a_{n}^{x}(s) 1_{[0,s]}^{\otimes n}(\cdot )1_{[a,T]}(s) \right) ^{\sim}(t_{1},
\ldots , t_{n+1})= \frac{1}{n+1} \sum _{i=1}^{n+1} a_{n}^{x}
(t_{i}) 1_{[0, t_{i}]} ^{\otimes n} (t_{1}, \ldots , \hat{t_{i}},
\ldots , t_{n+1})1_{[a,T]}(t_{i})
\end{equation*}
where as above $\hat{t_{i}}$ means that the variable $t_{i}$ is
missing.

\vskip0.2cm

We recall that $F$ is a random variable having the chaotic
decomposition $F=\sum_{n}I^{i}_{n}(f_{n})$ then its Sobolev-Watanabe
norm is given by
\begin{equation}
\label{wata} \Vert F\Vert ^{2}_{2,\alpha }=\sum_{n} (n +1)^{\alpha }
\Vert I_{n}(f_{n}) \Vert ^{2} _{L^{2} (\Omega)}.
\end{equation}

We will get
\begin{eqnarray*}
\Vert \xi(x) \Vert ^{2}_{2,\alpha }&=&\sum_{n\geq 0} (n+2)^{\alpha
} (n+1)! \Vert \left( a_{n}^{x}(s) 1_{[0,s]}^{\otimes n}(\cdot )
\right) ^{\sim} \Vert ^{2} _{L^{2}[0,T]^{n+1}} \\
&=&\sum_{n_{i}\geq 0} (n+2)^{\alpha } n! \int_{a}^{T} \left( p_{s}(x)
H_{n} \left( \frac{x}{\sqrt{s}}\right)\right) ^{2} ds .
\end{eqnarray*}
We use  the identity
\begin{equation}
\label{id} H_{n}(y)e^{-\frac{y^{2}}{2}} = (-1)^{[n/2]}
2^{\frac{n}{2}} \frac{2}{n!\pi } \int _{0}^{\infty } u^{n}
e^{-u^{2}} g(uy\sqrt{2}) du
\end{equation}
where $g(r)=\cos (r)$ if $n$ is even and $g(r)=\sin (r)$ if $n$ is
odd. Since $\vert g(r)\vert \leq 1$, we have the bound

\begin{equation}\label{cn}
\left|H_{n}(y)e^{-\frac{y^{2}}{2}} \right| \leq
2^{\frac{n}{2}}\frac{2}{n!\pi } \Gamma (\frac{n+1}{2}):=c_{n}
\end{equation}
Then
\begin{equation}
\label{b1} \Vert \xi(x) \Vert ^{2}_{2,\alpha } \leq cst. \sum
_{n}(n+2)^{\alpha }n! c_{n}^{2} \int _{a}^{T} \frac{1}{s}ds
\end{equation}
and this is finite for $\alpha <\frac{-1}{2}$ since by Stirling's
formula $n!c_{n}^{2}$ behaves as $cst. \frac{1}{\sqrt{n}}$.

We summarize the above discussion.

\begin{prop}
For any $x\in \mathbb{R}$ the functional $\xi (x)= \int _{0}^{T}
\delta (x-B_{s} ) dB_{s}$ belongs to the Sobolev -Watanabe space
$\mathbb{D}^{ -\alpha, 2 }$ for any $\alpha >\frac{1}{2}$.
\end{prop}

\begin{remark}
 As for the delta Dirac function, the regularity $\xi$ is the same with respect to $x$ and with respect to $\omega$.

\end{remark}

\subsection{The fractional case}
In this paragraph the driving process is a fractional Brownian
motion $(B^{H}_{t} )_{t\in [0,T]} $ with Hurst parameter $H\in
(\frac{1}{2},1)$.  We are interested to study the regularity as a functional
in Sobolev-Watanabe spaces of
\begin{equation}
\label{xiH} \xi (x) = \int_{0}^{T}  \delta \left( x-B^{H}_{s}\right)
dB^{H}_{s}
\end{equation}
when $x\in \mathbb{R}$ is fixed. We will that now the order of regularity in the Watanabe spaces changes and it differs from the order of regularity of the same functional with respect to the variable $x$.

We prove the following result.
\begin{prop}
Let $(B^{H}_{t} )_{t\in [0,T]} $ be a fractional Brownian motion with Hurst parameter $H\in
(\frac{1}{2},1)$. For any $x\in \mathbb{R}$ the functional $\xi (x)$ (\ref{xi}) is an element of the Sobolev -Watanabe space $\mathbb{D}^{-\alpha ,2}$ with $\alpha >\frac{3}{2}-\frac{1}{2H}$.
\end{prop}
{\bf Proof: } We will have in this case
\begin{equation*}
\Vert \xi(x) \Vert ^{2} _{2, \alpha } = \sum _{n\geq 0} (n+2)
^{\alpha } (n+1) ! \Vert \left( a_{n}^{x}(s)1_{[0,s]}^{\otimes
n}(\cdot )\right) ^{\sim}  \Vert ^{2} _{{\cal{H}}_{H}^{\otimes n+1}}
\end{equation*}
where we denoted  by ${\cal{H}}_{H}$ the canonical Hilbert space of the
fractional Brownian motion.

As in the proof of Proposition 2 we obtain
\begin{eqnarray*}
&&\Vert \left( a_{n}^{x}(s)1_{[0,s]}^{\otimes
n}(\cdot )\right) ^{\sim}  \Vert ^{2} _{{\cal{H}}_{H}^{\otimes n+1}} \\
&=& \frac{ (H(2H-1))
^{n+1}}{(n+1)^{2}} \sum _{i,j=1} ^{n+1} \int _{[0,T]^{n+1}} \int
_{[0,T]^{n+1}}du_{1}\ldots du_{n+1}dv_{1}\ldots dv_{n+1} \prod _{l=1}^{n+1} \vert u_{l}-v_{l} \vert ^{2H-2} \\
&&\times a_{n}^{x}(u_{i} ) a_{n}^{x} (v_{j} )1_{[0,u_{i}]}^{\otimes
n}(u_{1}, \ldots , \hat{u}_{i}, \ldots ,
u_{n+1})1_{[0,v_{j}]}^{\otimes n}(v_{1}, \ldots , \hat{v}_{j},
\ldots , v_{n+1})\\
&=&\frac{ (H(2H-1)) ^{n+1}}{(n+1)^{2}} \sum_{i=1} ^{n+1}\int
_{[0,T]^{n+1}} \int
_{[0,T]^{n+1}}du_{1}\ldots du_{n+1} dv_{1}\ldots dv_{n+1}\prod _{l=1}^{n+1} \vert u_{l}-v_{l} \vert ^{2H-2} \\
&&\times  a_{n}^{x}(u_{i} ) a_{n}^{x} (v_{i} )1_{[0,u_{i}]}^{\otimes
n}(u_{1}, \ldots , \hat{u}_{i}, \ldots ,
u_{n+1})1_{[0,v_{i}]}^{\otimes n}(v_{1}, \ldots , \hat{v}_{i},
\ldots , v_{n+1})\\
&&+\frac{ (H(2H-1)) ^{n+1}}{(n+1)^{2}} \sum _{i\not=j; i,j=1} ^{n+1}
\int _{[0,T]^{n+1}} \int
_{[0,T]^{n+1}}du_{1}\ldots du_{n+1}dv_{1}\ldots dv_{n+1} \prod _{l=1}^{n+1} \vert u_{l}-v_{l} \vert ^{2H-2} \\
&&\times a_{n}^{x}(u_{i} ) a_{n}^{x} (v_{j} )1_{[0,u_{i}]}^{\otimes
n}(u_{1}, \ldots , \hat{u}_{i}, \ldots ,
u_{n+1})1_{[0,v_{j}]}^{\otimes n}(v_{1}, \ldots , \hat{v}_{j},
\ldots , v_{n+1})\\
&:=& A(n)+B(n).
\end{eqnarray*}
The first term $A(n)$ equals, by symmetry,
\begin{equation*}
A(n)=\frac{ (H(2H-1)) ^{n+1}}{n+1}\int_{0}^{T} \int_{0}^{T} dudv
a_{n}^{x}(u) a_{n}^{x}(v)  \vert u-v\vert ^{2H-2} \left(   \int _{0}
^{u}\int _{0}^{v} \vert u' -v' \vert ^{2H-2}  du'dv' \right) ^{n}
\end{equation*}
and using equality (\ref{II}) we get
\begin{eqnarray*}
A(n)&=&\frac{ (H(2H-1)) }{n+1}\int_{0}^{T}
\int_{0}^{T}dudv\\
&&\times \vert u-v\vert ^{2H-2} R(u,v) ^{n} u^{-Hn} v^{-Hn}
p_{u^{2H}}(x) p_{v^{2H}}(x)H_{n} \left( \frac{x}{u^{H}}\right)H_{n}
\left( \frac{x}{v^{H}}\right).
\end{eqnarray*}
 Using the identity (\ref{id})
\begin{eqnarray*}
&&\sum_{n\geq 0} (n+2) ^{\alpha } (n+1)! A(n) \\
&=& c(H) \sum _{n\geq 0} (n+2)^{\alpha } n!\int_{0}^{T}
\int_{0}^{T}dudv\vert u-v\vert ^{2H-2}\\
&&\times R(u,v) ^{n} u^{-Hn} v^{-Hn} p_{u^{2H}}(x)
p_{v^{2H}}(x)H_{n} \left( \frac{x}{u^{H}}\right)H_{n} \left(
\frac{x}{v^{H}}\right)\\
&=&c(H) \sum _{n\geq 0} (n+2)^{\alpha }
n!c_{n}^{2}\int_{0}^{T}
\int_{0}^{T}dudv\vert u-v\vert ^{2H-2}R(u,v) ^{n} u^{-Hn} v^{-Hn} u^{-H}v^{-H}.
\end{eqnarray*}
By the selfsimilarity of the fBm we have $R(u,v)= u^{-2H} R(1,
\frac{v}{u})$ and by the change of variables  $v=zu $ we obtain

\begin{eqnarray}
&&\sum_{n\geq 0} (n+2) ^{\alpha } (n+1)! A(n)
\leq \nonumber\\
&& c(H)\sum _{n\geq 0} (n+2) ^{\alpha } n!c_{n}^{2} \left( \int _{0}^{T}
u^{1-2H}du\right) \int _{0}^{1} \frac{ R(1,z) ^{n}}{z^{Hn}}
\frac{(1-z)^{2H-2}}{z^{H}}dz .\label{R11}
\end{eqnarray}
Using Lemma 2 (actually a slightly modification of it) in \cite{Ed}
\begin{equation}
\label{edd} \int _{0}^{1} \frac{ R(1,z) ^{n}}{z^{Hn}}
\frac{(1-z)^{2H-2}}{z^{H}}dz \leq c(H) \frac{1}{n^{\frac{1}{2H}}}
\end{equation}
and the right hand side of (\ref{R11}) is bounded, modulo a constant, by
$$\sum _{n\geq 0} (n+2) ^{\alpha }n!c_{n} ^{2} n^{-\frac{1}{2H}}$$ and since
$n! c_{n}^{2}$
behaves as $\frac{1}{\sqrt{n}}$, the last sum is convergent if
$\frac{1}{2H} -\alpha +\frac{1}{2} >1$, or $-\alpha
>-\frac{1}{2H}+\frac{1}{2} $.

\vskip0.2cm

 Let us regard now the sum involving the term $B(n)$. It will actually decide the regularity of the functional $\xi(x)$.
Following the computations contained in the proof of Proposition \ref{p3}
\begin{eqnarray*}
&& \sum _{n\geq 0}(n+2) ^{\alpha } (n+1)! B(n) \\
&=&c(H) \sum _{n\geq 0}(n+2) ^{\alpha } (n+1)!  \frac{n(n+1)}{(n+1)
^{2}} \int _{0}^{T}\int _{0}^{T}\int _{0}^{T}\int
_{0}^{T}du_{1}du_{2} dv_{1}dv_{2} \\
&& \times 1_{[0,u_{1}]} (u_{2})1_{[0,v_{2}]} (v_{1})R(u_{1}, v_{2})
^{n-1} \vert u_{1}-v_{1}\vert ^{2H-2}\vert u_{2}-v_{2}\vert
^{2H-2}a_{n}^{x}(u_{1}) a_{n}^{x}(v_{2})\\
&\leq &c(H) \sum _{n\geq 0}(n+2) ^{\alpha } n n! c_{n}^{2}  \\
&&\times \int _{0}^{T}\int _{0}^{T} du_{1}dv_{2}R(u_{1}, v_{2})
^{n-1} u_{1}^{-Hn}v_{2}^{-Hn}  \left( \int_{0}^{u_{1}}\vert
u_{2}-v_{2}\vert ^{2H-2} du_{2}\right) \left( \int _{0}^{v_{2}}\vert
u_{1}-v_{1}\vert ^{2H-2} \right)\\
&\leq &c(H) \sum _{n\geq 0}(n+2) ^{\alpha } n n! c_{n}^{2}  \int _{0}^{T}\int _{0}^{T}
du_{1}dv_{2}R(u_{1}, v_{2}) ^{n-1} u_{1}^{-Hn}v_{2}^{-Hn}.
\end{eqnarray*}
Using again Lemma 2 in \cite{Ed} , we get that the integral
\begin{equation*}
\int _{0}^{T}\int _{0}^{T} du_{1}dv_{2}R(u_{1}, v_{2}) ^{n-1}
u_{1}^{-Hn}v_{2}^{-Hn} \leq c(H) n^{-\frac{1}{2H} }
\end{equation*}
and since the sequence $n!c_{n}^{2}$ behaves when $n\to \infty $ as $\frac{1}{\sqrt{n}}$
we obtain that the sum $\sum _{n\geq 0}(n+2) ^{\alpha } (n+1)! B(n)$
converges if $- \alpha - \frac{1}{2}+ \frac{1}{2H}>1$ and this gives
$-\alpha >\frac{3}{2}-\frac{1}{2H}$. \qed

\begin{remark}
Only when $H=\frac{1}{2}$ we retrieve the same order of regularity of (\ref{xiH}) as  a function of $x$ and as a function of $\omega$.

\end{remark}




\vskip0.5cm




\end{document}